%




\font\sc=cmcsc10 \rm


\newcount\secnb 
\newtoks\secref 
\secnb=0\secref={}
\newcount\subnb 
\newtoks\subref 
\subnb=0\subref={}
\newcount\parnb 
\parnb=0

\def\smallskip{\par\vskip 2mm}
\def\medskip{\par\vskip 5mm}
\def\goodbreak{\penalty -100}

\def\introduction{\subnb=0\parnb=0
\medskip\goodbreak
\centerline{\bf Introduction}
\smallskip\nobreak}

\def\section#1{\global\advance\secnb by 1
\secref=\expandafter{\the\secnb}
\subref={}
\subnb=0\parnb=0
\vfill\eject
\centerline{\bf\S\the\secref.\ #1}
\smallskip\nobreak}

\def\subsection#1{\global\advance\subnb by 1
\secref=\expandafter{\the\secnb.}
\subref=\expandafter{\the\subnb}
\parnb=0
\smallskip\goodbreak
\leftline{\bf\the\secref\the\subref.\ #1}
\par\nobreak}

\def\references{\vfill\eject
\centerline{\bf References}
\smallskip\nobreak}

\def\paragraph#1{\global\advance\parnb by 1
\smallskip\goodbreak\noindent\the\secref\the\subref.\the\parnb. {#1}\par\nobreak}

\def\proof{\smallskip\goodbreak\noindent{\it Proof:}\par\nobreak}
\def\endproof{\relax}


\def\label#1{\relax}

\def\ref#1{\csname crossref#1\endcsname}


\newcount\bibrefnb
\bibrefnb=0

\def\bibitem#1{\global\advance\bibrefnb by 1\item{\the\bibrefnb.}}

\def\refto#1{\setbox1=\hbox{\sc #1}\copy1}
\def\reftosame{\hbox to \the\wd1{\hrulefill}}

\def\cite#1{\csname bibref#1\endcsname}

\def\bye{\medskip\leftline{{\it Mail address:}}{\leftskip=1cm\parindent=0cm

Laboratoire J.A. Dieudonn\'e,\hfill\break
Universit\'e de Nice,\hfill\break
Parc Valrose,\hfill\break
F-06108 Nice Cedex 02 (France).\par}

\smallskip\leftline{{\it E-mail address:}}{\leftskip=1cm\parindent=0cm

Clemens Berger $<$cberger@math.unice.fr$>$,\hfill\break
Benoit Fresse $<$fresse@math.unice.fr$>$.\par}

\vfill\supereject\end}


\font\tenms=msbm10
\font\sevenms=msbm7
\font\fivems=msbm5

\newfam\msfam
\textfont\msfam=\tenms
\scriptfont\msfam=\sevenms
\scriptscriptfont\msfam=\fivems
\def\ms{\fam\msfam}

\def\N{{\ms N}}
\def\Z{{\ms Z}}

\def\F{{\ms F}}

\def\pt{\mathop{pt}}

\def\dg{\mathop{dg}}

\def\Set{{{\cal S}et}}
\def\sSet{{\cal S}}

\def\Mod{\mathop{{\cal M}od}\nolimits}

\def\Hom{\mathop{\rm Hom}\nolimits}

\def\lim{\mathop{\rm lim}}

\def\hocolim{\mathop{\rm colim}}


\def\ker{\mathop{\rm ker}\nolimits}

\def\rto#1#2{\smash{\mathop{\hbox to 6 mm{\rightarrowfill}}
\limits^{\scriptstyle #1}_{\scriptstyle #2}}}
\def\lto#1#2{\smash{\mathop{\hbox to 6 mm{\leftarrowfill}}
\limits^{\scriptstyle #1}_{\scriptstyle #2}}}

\def\O{{\cal O}}
\def\P{{\cal P}}
\def\Q{{\cal Q}}
\def\A{{\cal A}} 
\def\C{{\cal C}} 
\def\D{{\cal D}} 
\def\E{{\cal E}} 
\def\X{{\cal X}} 
\def\Z{{\cal Z}} 
\def\W{{\cal W}} 

\def\AW{\mathop{AW}}
\def\EZ{\mathop{EZ}}
\def\TR{\mathop{TR}}
\def\TC{\mathop{TC}}
\def\id{\mathop{id}}

\def\seg{\hbox to 6mm{\hrulefill}}
\def\length{\mathop{\rm length}\nolimits}

\def\sgn{\mathop{sgn}}

\input xy
\xyoption{all}

%
%

\magnification=1200

\centerline{\bf Combinatorial operad actions on cochains}
\smallskip
\centerline{Clemens Berger}\centerline{Benoit Fresse}
\centerline{October 8, 2002}
\medskip

\medskip\centerline{\bf Abstract}

\smallskip{\lineskiplimit=2pt\lineskip=2pt\baselineskip=10pt\leftskip=1cm\rightskip=1cm\sevenrm
A classical E-infinity operad is formed by the bar construction of the symmetric groups.
Such an operad has been introduced by M. Barratt and P. Eccles
in the context of simplicial sets in order to have an analogue of the Milnor FK-construction
for infinite loop spaces.
The purpose of this article is to prove that the associative algebra structure
on the normalized cochain complex of a simplicial set
extends to the structure of an algebra over the Barratt-Eccles operad.
We also prove that differential graded algebras over the Barratt-Eccles operad form a closed model category.
Similar results hold for the normalized Hochschild cochain complex of an associative algebra.
More precisely, the Hochschild cochain complex is acted on by a suboperad of the Barratt-Eccles operad
which is equivalent to the classical little squares operad.\par}

\introduction

\smallskip
The main purpose of this article is to provide the normalized cochain complex of a simplicial set $N^*(X)$
with the structure of an algebra over a suitable $E_\infty$-operad.
Our construction is motivated by results of M. Mandell ({\it cf}. [\cite{M}])
which imply that any such structure determines the homotopy type of $X$,
provided $X$ is nilpotent $p$-complete and has a finite $p$-type.

To be precise,
we consider a classical $E_\infty$-operad called the Barratt-Eccles operad,
because it comes from a simplicial operad introduced by M. Barratt and P. Eccles
for the study of infinite loop spaces ({\it cf}. [\cite{BE}]).
This operad, which we denote by $\E$, consists of a sequence of chain complexes $\E(r)$, $r\in\N$,
formed by the normalized bar constructions
of the symmetric groups $\Sigma_r$.
An $\E$-algebra is an $\F$-module $A$
equipped with evaluation products $\E(r)\otimes A^{\otimes r}\,\rightarrow\,A$
defined for $r\in\N$.
Accordingly,
an element $\rho\in\E(r)$ gives rise to an operation
$\rho: A^{\otimes r}\,\rightarrow\,A$.

We obtain the following result:

\smallskip\noindent{\sc Theorem}\par\nobreak
{\it For any simplicial set $X$,
we have evaluation products
$\E(r)\otimes N^*(X)^{\otimes r}\,\rightarrow\,N^*(X)$
(functorial in $X$)
which give the normalized cochain complex $N^*(X)$
the structure of an algebra over the Barratt-Eccles operad $\E$.
In particular, the classical cup-product of cochains
is an operation $\mu_0: N^*(X)^{\otimes 2}\,\rightarrow\,N^*(X)$
associated to an element $\mu_0\in\E(2)_0$.}

\smallskip
The existence of such $\E$-algebra structures is proved by Justin Smith in [\cite{JuS}].
In this article,
we follow a different and more combinatorial approach in order to make the evaluation products
$\E(r)\otimes N^*(X)^{\otimes r}\,\rightarrow\,N^*(X)$
completely explicit.
As an example,
in the case of spheres $S^n$,
we obtain a simple formula for the evaluation products
$$\E(r)\otimes N^*(S^n)^{\otimes r}\,\rightarrow\,N^*(S^n)$$
({\it cf}. theorem \ref{SphereModels}).

The Barratt-Eccles operad comes equipped with a diagonal $\E(r)\,\rightarrow\,\E(r)\otimes\E(r)$
because it is the chain complex of a simplicial operad.
This structure has interesting applications
which make the Barratt-Eccles operad suitable for calculations in homotopy theory.
For instance,
we have a straighforward construction of path-objects
in the category of $\E$-algebras.
As a result,
we prove that differential graded algebras over the Barratt-Eccles operad
form a closed model category.
We do not know whether this property is satisfied for any $E_\infty$-operad.

Here is a general result obtained in section 3.1.
The Barratt-Eccles operad $\P = \E$ is a $\Sigma_*$-projective resolution
of the operad $\Q = \C$ of commutative algebras.
In general,
we may replace an operad $\Q$ by the operad tensor product $\P = \E\otimes\Q$
which form a $\Sigma_*$-projective resolution of $\Q$.
We obtain the following theorem:

\smallskip\noindent{\sc Theorem}\par\nobreak
{\it We consider an operad $\P$ such that $\P = \E\otimes\Q$.
We assume the classical definition of a weak equivalence (respectively, of a fibration)
in the category of $\P$-algebras.
Namely,
a weak equivalence (respectively, a fibration)
is an algebra morphism which is a weak equivalence (respectively, a fibration)
in the category of dg-modules.
We claim that these definitions provide the category of $\P$-algebras
with the structure of a closed model category.}

\smallskip
The simplicial Barratt-Eccles operad has a filtration
$$F_1\W\subseteq F_2\W\subseteq\cdots\subseteq F_n\W\subseteq\cdots\subseteq\W$$
by simplicial suboperads $F_n\W$
whose topological realizations $|F_n\W|$ are equivalent
to the classical operads of little $n$-cubes
({\it cf}. [\cite{Bg}], [\cite{K}], [\cite{JfS}]).
We have an induced filtration on the associated differential graded operad
$$F_1\E\subseteq F_2\E\subseteq\cdots\subseteq F_n\E\subseteq\cdots\subseteq\E.$$
The operad $F_1\E$ is identified with the associative operad.
The operad $F_2\E$ is the chain operad of a simplicial operad equivalent to the little squares operad.
The next result implies that this operad $F_2\E$ gives a solution
to a conjecture of Deligne:

\smallskip\noindent{\sc Theorem}\par\nobreak
{\it For any associative algebra $A$,
the normalized Hochschild cochain complex of $A$ (denoted by $C^*(A,A)$)
is endowed with natural evaluation products
$F_2\E(r)\otimes C^*(A,A)^{\otimes r}\,\rightarrow\,C^*(A,A)$
which give $C^*(A,A)$ the structure of an algebra over the operad $F_2\E$.
In particular, the classical cup-product of Hochschild cochains is an operation
$\mu_0: C^*(A,A)^{\otimes 2}\,\rightarrow\,C^*(A,A)$
associated to an element $\mu_0\in\E(2)_0$.}

\smallskip
The pivot of our constructions is formed by a quotient operad $\X$
of the Barratt-Eccles operad $\E$
which we call the surjection operad.
This operad is also introduced (with different sign conventions) by J. McClure and J. Smith
in their work on the Deligne conjecture ({\it cf}. [\cite{MCSPrism}], [\cite{MCSSeq}])
to which we refer in the article.
To be precise,
the Barratt-Eccles operad $\E$ operates on cochain complexes of simplicial sets
and on Hochschild complexes of associative algebras
through a morphism $\TR: \E\,\rightarrow\,\X$
which we call the table reduction morphism.

One observes precisely that the surjection operad $\X$ acts naturally on the cochain complex of a simplicial set.
This idea goes back to D. Benson ({\it cf}. [\cite{Bn}, Section 4.5])
and to the original construction of the reduced square operations
by N. Steenrod ({\it cf}. [\cite{S}]).
The introduction of the surjection operad in the context of the Deligne conjecture
is due to J. McClure and J. Smith.
According to these authors,
the surjection operad $\X$ is endowed with a filtration
$$F_1\X\subseteq F_2\X\subseteq\cdots\subseteq F_n\X\subseteq\cdots\subseteq\X$$
by suboperads $F_n\X$ such that $F_1\X$ is identified with the associative operad.
In addition,
the operad $F_2\X$ acts naturally on the Hochschild cochain complex of an associative algebra
({\it cf}. [\cite{MCSPrism}], [\cite{MCSSeq}]).

Finally,
we deduce our main theorems from the following result:

\smallskip\noindent{\sc Theorem}\par\nobreak
{\it The surjection operad $\X$ is a quotient of the Barratt-Eccles operad $\E$.
More precisely,
we have a surjective morphism of filtered operads $\TR: \E\,\rightarrow\,\X$
such that $F_1\TR: F_1\E\,\rightarrow\,F_1\X$
is the identity isomorphism of the associative operad.}

\smallskip
We give a geometric interpretation of the table reduction morphism 
$\TR: \E\,\rightarrow\,\X$ in the note [\cite{BF}].

J. McClure and J. Smith prove that the surjection operad $F_n\X$
is equivalent to the chain-operad of little $n$-cubes.
(As a consequence, the operad $F_2\X$ is also a solution to the Deligne conjecture.)
We prove more precisely that the table reduction morphism induces weak equivalences
$F_n\TR: F_n\E\,\rto{\sim}{}\,F_n\X$.
But, since the operads $F_n\X$ are not directly associated to a simplicial operad,
it is not clear that certain of our constructions are still available
in the category of algebras over the surjection operad.
This motivates the introduction of the Barratt-Eccles operad
and of the table reduction morphism.

\smallskip
We refer to section 0 for conventions on operads
which are adopted throughout the article.
The detailed definition of the Barratt-Eccles operad is given in section 1.1.
The surjection operad is introduced in section 1.2.
Section 1.3 is devoted to the construction of the table reduction morphism.
The properties of the morphism are proved in sections 1.4 and 1.5.
(These two sections may be skipped in a first reading.)
The little cubes filtrations and the applications to the Deligne conjecture are explained in section 1.6.
The construction of the algebra structure on the normalized cochain complex of a simplicial set
is achieved in section 2.
The applications of our constructions to homotopical algebra are explained in section 3.

\medskip\centerline{\bf Contents}

\smallskip{\parindent=0cm\obeylines

0. Conventions

1. The operads

{\parindent=1cm
1.1. The Barratt-Eccles operad
1.2. The surjection operad
1.3. The table reduction morphism
1.4. The operad morphism properties
1.5. On composition products
1.6. On little cubes filtrations}

2. The interval cut operations on chains

{\parindent=1cm
2.1. The Eilenberg-Zilber operad
2.2. The construction of interval cut operations}

3. On closed model structures

{\parindent=1cm
3.1. The closed model structure
3.2. On spheres, cones and suspensions
3.3. Some proofs}

}

\medskip\noindent{\it Acknowledgements:}
The morphisms from the Barratt-Eccles operad to the surjection operad and to the Eilenberg-Zilber operad
have been constructed (modulo the signs) by the second author
and lectured on at the fall school ``Lola'' in Malaga in November 2000.
The second author would like to thank the organizers and the participants of the school
who provided the first motivation for this work.
We thank Jim McClure and Jeff Smith for their interest on our work
which, at some places, should complement their results on the Deligne conjecture.
In fact, the first author realized that the Barratt-Eccles operad provides a solution to the Deligne conjecture
after reading the announcement of the article [\cite{MCSSeq}].

\secnb=-1\section{Conventions}

\paragraph{\it On permutations}
We let $\Sigma_r$, $r\in\N$, denote the permutation groups.
The element $1_r\in\Sigma_r$ is the identity permutation $1_r = (1,2,\ldots,r)$.
In general, a permutation $\sigma\in\Sigma_r$ is specified
by the sequence $(\sigma(1),\sigma(2),\ldots,\sigma(r))$
formed by its values.

Given $s_1,\ldots,s_r\in\N$,
we have a block permutation $\sigma_*(s_1,\ldots,s_r)\in\Sigma_{s_1+\cdots+s_r}$
associated to any $\sigma\in\Sigma_r$
and which is determined by the following process:
in the sequence $(\sigma(1),\sigma(2),\ldots,\sigma(r))$,
we replace the occurence of each value $k = 1,\ldots,r$
by the sequence
$$s_1+s_2+\cdots+s_{k-1}+1,s_1+s_2+\cdots+s_{k-1}+2,\ldots,s_1+s_2+\cdots+s_{k-1}+s_k.$$
For instance,
if $\sigma$ is the transposition $\sigma = (2,1)$,
then $\sigma_*(p,q)$ is the block transposition
$\sigma_*(p,q) = (p+1,\ldots,p+q,1,\ldots,p)$.

\paragraph{\it On the symmetric monoidal category of dg-modules}
Let us fix a ground ring $\F$.
The category of $\F$-modules is denoted by $\Mod_\F$.
The category of differential graded modules is denoted by $\dg\Mod_\F$.

If $V$ is a differential graded module (or, for short, a dg-module),
then $V$ has either an upper grading $V = V^*$ or a lower $V= V_*$ grading.
The equivalence between a lower and an upper grading is given by the classical rule.
We have explicitly $V_{d} = V^{-d}$.
Let us mention that differential degrees are assumed to range over the integers
without any restriction.
In general,
the differential of a dg-module is denoted by $\delta: V^*\,\rightarrow\,V^{*+1}$.

We equip the category of dg-modules with the classical tensor product of dg-modules
$$\otimes: \dg\Mod_\F\times\dg\Mod_\F\,\rightarrow\,\dg\Mod_\F.$$
To be explicit, if $U,V\in\dg\Mod_\F$, then $U\otimes V\in\dg\Mod_\F$ has the module
$$(U\otimes V)^n = \bigoplus_{d\in{\ms Z}} U^d\otimes V^{n-d}$$
in degree $n$.
The differential
$\delta: U\otimes V\,\rightarrow\,U\otimes V$
is given by the classical derivation rule
$$\delta(u\otimes v) = \delta(u)\otimes v + (-1)^d u\otimes\delta(v),$$
for all homogeneous elements $u\in U^d$ and $v\in V^e$.
The symmetry operator $\tau: U\otimes V\,\rightarrow\,V\otimes U$ follows the sign-rule.
Precisely,
a permutation of homogeneous symbols,
which have respective degrees $d$ and $e$,
produces a sign $(-1)^{d e}$.
In particular,
given homogeneous elements $u\in U^d$ and $v\in V^e$,
we have
$$\tau(u\otimes v) = (-1)^{d e} v\otimes u.$$
Similarly,
in the derivation rule above,
the permutation of the differential $\delta$, which has degree $1$,
with the element $u$, homogeneous of degree $d$,
produces the sign $(-1)^d$.

A morphism $f: U\,\rightarrow\,V$
is homogeneous of degree $d$ if $f(U^*)\subset V^{*+d}$,
for all $*\in{\ms Z}$.
The differential $\delta(f): U\,\rightarrow\,V$ is given by the commutator of $f$
with the differentials of $U$ and $V$:
$$\delta(f) = \delta\cdot f - (-1)^d f\cdot\delta.$$
The dg-module formed by the homogeneous morphisms is denoted by $\Hom^*_\F(U,V)$.
Thus, we have explicitly:
$$\Hom^d_\F(U,V)
= \prod_{*\in{\ms Z}}\Hom_\F(U^{*},V^{*+d})
= \prod_{d\in{\ms Z}}\Hom_\F(U_{*+d},V_{*}).$$
This dg-module $\Hom^*_\F(U,V)$ is an internal hom
in the category of dg-modules.

\paragraph{\it The symmetric monoidal categories}
We are concerned with the following classical symmetric monoidal categories:
\item{1)} the category of $\F$-modules $\Mod_\F$
equipped with the classical tensor product
$$\otimes: \Mod_\F\times\Mod_\F\,\rightarrow\,\Mod_\F;$$
\item{2)} the category of dg-modules $\dg\Mod_\F$
equipped with the tensor product of dg-modules
$$\otimes: \dg\Mod_\F\times\dg\Mod_\F\,\rightarrow\,\dg\Mod_\F,$$
whose definition is recalled in the paragraph above;
\item{3)} the category of sets $\Set$ equipped with the cartesian product
$$\times: \Set\times\Set\,\rightarrow\,\Set;$$
\item{4)} the category of simplicial sets $\sSet$
equipped with the cartesian product of simplicial sets
$$\times: \sSet\times\sSet\,\rightarrow\,\sSet.$$

In addition,
we shall consider the symmetric monoidal category of (differential graded) coalgebras.
Classically,
a coalgebra is a (differential graded) module $X\in\Mod_\F$
together with an associative diagonal
$\Delta: X\,\rightarrow\,X\otimes X$.
One observes that a tensor product of coalgebras $X\otimes Y\in\Mod_\F$
is equipped with a canonical diagonal
which is given by the composite arrow
$$X\otimes Y
\,\rto{\Delta\otimes\Delta}{}\,X\otimes X\otimes Y\otimes Y
\,\rto{(1,3,2,4)}{}\,X\otimes Y\otimes X\otimes Y.$$
Therefore, the dg-coalgebras form a symmetric monoidal category.

\paragraph{\it On normalized chains and cochains}\label{NormChains}
If $X$ is a simplicial set,
then $C_*(X)$ denotes the dg-module
whose degree $d$ component $C_d(X)$ is the free $\F$-module
generated by the set $X_d$
of the $d$-dimensional simplices of $X$.
The differential of $x\in X_d$ in $C_*(X)$ is given by the classical formula:
$\delta(x) = \sum_{i=0}^{d} (-1)^i d_i(x)$.
The normalized chain complex $N_*(X)$ is the quotient of the dg-module $C_*(X)$
by the degeneracies:
$$N_d(X) = C_d(X)/s_0 C_{d-1}(X)+\cdots+s_{d-1} C_{d-1}(X).$$
We consider also the dual cochain complexes $C^*(X) = \Hom_\F(C_*(X),\F)$
and $N^*(X) = \Hom_\F(N_*(X),\F)$.

Let us recall that the normalized chain complex $X\mapsto N_*(X)$
is a (symmetric) monoidal functor from the category of simplicial sets
to the category of dg-modules.
More precisely,
the Eilenberg-Zilber morphism provides a natural equivalence
$$\EZ:\,N_*(X)\otimes N_*(Y)\,\rto{\sim}{}\,N_*(X\times Y)$$
which is associative and commutes with the symmetry isomorphism.

In addition,
the dg-module $N_*(X)$ is equipped with a natural diagonal
(the Alexander-Whitney diagonal)
and forms a dg-coalgebra.
We conclude that the normalized chain complex $X\mapsto N_*(X)$
is a (symmetric) monoidal functor from the category of simplicial sets
to the category of dg-coalgebras.

We refer to the classical textbook of S. Mac Lane ({\it cf}. [\cite{ML}, Section VIII.8])
for the properties of the Eilenberg-Zilber equivalence.

\paragraph{\it On the notion of an operad}
The main structures involved in our article are operads
in the symmetric monoidal category of dg-modules.
But, for simplicity,
we recall the definition of an operad in the category of $\F$-modules.
In fact,
the notion of an operad makes sense in any symmetric monoidal category.
For an introduction to this subject,
we refer to the articles of E. Getzler and J. Jones ({\it cf}. [\cite{GeJ}]),
V. Ginzburg and M. Kapranov ({\it cf}. [\cite{GiK}])
and to the work of P. May ({\it cf}. [\cite{MLoop}]).

A $\Sigma_*$-module $M$ is a sequence $M = M(r)$, $r\in\N$,
where $M(r)$ is a representation of the permutation group $\Sigma_r$.
A morphism of $\Sigma_*$-modules $f: M\,\rightarrow\,N$
is a sequence of representation morphisms
$f: M(r)\,\rightarrow\,N(r)$.

An operad $\P$ is a $\Sigma_*$-module equipped with composition products
$$\P(r)\otimes\P(s_1)\otimes\cdots\otimes\P(s_r)\,\rightarrow\,\P(s_1+\cdots+s_r)$$
defined for $r\geq 1$, $s_1,\ldots,s_r\geq 0$,
and which statisfies some natural equivariance and associativity properties.
A morphism of operads $f: \P\,\rightarrow\,\Q$ is a morphism of $\Sigma_*$-modules
which commutes with the operad composition products.

In general,
an element $p\in\P(r)$ represents a (multilinear) operation
in $r$ variables $x_1,\ldots,x_r$.
Thus, we may write $p = p(x_1,\ldots,x_r)$.
The action of a permutation $\sigma\in\Sigma_r$ on $p\in\P(r)$ is given by the permutation of variables.
Explicitly,
we have $\sigma\cdot p = p(x_{\sigma(1)},\ldots,x_{\sigma(r)})$.
The composition product of $p\in\P(r)$ and $q_1\in\P(s_1),\ldots,q_r\in\P(s_r)$
is denoted by $p(q_1,\ldots,q_r)\in\P(s_1+\cdots+s_r)$.
In fact,
this operation is given by the substitution of variables $x_1,\ldots,x_r$
by operations $q_1,\ldots,q_r$.

We assume that the composition product of an operad has a unit $1\in\P(1)$.
This unit represents the identity operation $1(x_1) = x_1$.
In this situtation,
there are partial composition products
$$\circ_k: \P(r)\otimes\P(s)\,\rightarrow\,\P(r+s-1)$$
defined for $k = 1,\ldots,r$
and which determine the operad composition products.
Explicitly,
the partial composition products of $p\in\P(r)$ and $q\in\P(s)$
are defined by the relation $p\circ_k q = p(1,\ldots,1,q,1,\ldots,1)$
(the operation $q$ is set at the $k$th entry of $p$).

\paragraph{\it On algebras over an operad}
The operad $\P$ has an associated category of algebras.
By definition,
a $\P$-algebra is an $\F$-module $A$
equipped with equivariant evaluation products
$$\P(r)\otimes A^{\otimes r}\,\rightarrow\,A$$
defined for $r\geq 0$.
This evaluation product has to be associative with respect
to the operad composition product
and unital with respect to the operad unit.
The evaluation product of an operation $p\in\P(r)$ at $a_1,\ldots,a_r\in A$
is denoted by $p(a_1,\ldots,a_r)\in A$.
The equivariance relation reads
$(\sigma\cdot p)(a_1,\ldots,a_r) = p(a_{\sigma(1)},\ldots,a_{\sigma(r)})$.
A morphism of $\P$-algebras $f: A\,\rightarrow\,B$ is a morphism of $\F$-modules
which commutes with the evaluation products.
The free $\P$-algebra generated by an $\F$-module $V$ is denoted by $\P(V)$.

\paragraph{\it The endomorphism operad}\label{EndOp}
The endomorphism operad of an $\F$-module $V$ is the operad
such that $\Hom_V(r) = \Hom_\F(V^{\otimes r},V)$
and whose composition product is defined by composing multilinear maps.
If $f: V^{\otimes r}\,\rightarrow\,V$, then the map $\sigma\cdot f: V^{\otimes r}\,\rightarrow\,V$
is defined by $f(v_1,\ldots,v_r) = f(v_{\sigma(1)},\ldots,v_{\sigma(r)})$.
An operad morphism $\P\,\rightarrow\,\Hom_V$
gives the dg-module $V$ the structure of a $\P$-algebra
({\it cf}. [\cite{GiK}]).
Accordingly,
the operad $\Hom_V$ is the universal operad
such that the dg-module $V$ has the structure of a $\Hom_V$-algebra.

\paragraph{\it On the normalized chain complex of a simplicial operad}
If $\O$ is a simplicial operad,
then the associated normalized chain complexes $\P(r) = N_*(\O(r))$
form an operad in the category of dg-coalgebras,
because the normalized chain complex functor $N_*(-)$ is symmetric monoidal
as reminded in the paragraph \ref{NormChains}.
Explicitly,
the Alexander-Whitney diagonal
$$N_*(\O(r))\,\rightarrow\,N_*(\O(r))\otimes N_*(\O(r))$$
gives $N_*(\O(r))$ the structure of a dg-coalgebra
and the composite of the composition product of $\O$
with the Eilenberg-Zilber equivalence
$$\eqalign{ N_*(\O(r))\otimes N_*(\O(s_1))\otimes\cdots\otimes N_*(\O(s_r))
& \,\rightarrow\,N_*(\O(r)\times\O(s_1)\times\cdots\times\O(s_r)) \cr
& \,\rightarrow\,N_*(\O(s_1+\cdots+s_r)) \cr }$$
gives the structure of a dg-operad.

\paragraph{\it On Hopf operads}
An operad $\P$ in the monoidal category of coalgebras is called a Hopf operad (as in [\cite{GeJ}]).
Accordingly, a dg-operad $\P$ is equipped with the structure of a Hopf operad
provided we have diagonals
$$\P(r)\,\rightarrow\,\P(r)\otimes\P(r)$$
that commute with the operad composition products
$$\P(r)\otimes\P(s_1)\otimes\cdots\otimes\P(s_r)\,\rightarrow\,\P(s_1+\cdots+s_r).$$
Similarly,
a $\P$-algebra together with a diagonal
is called a Hopf $\P$-algebra.

In this article, we consider only $\P$-algebras in the category of dg-modules (with no coalgebra structure).
We are just interested in the following property of algebras over a Hopf operad:
a tensor product of $\P$-algebras $A\otimes B\in\Mod_\F$
is equipped with a natural $\P$-algebra structure.
Precisely,
the composite
$$\eqalign{ \P(r)\otimes(A\otimes B)^{\otimes r} &
\,\rto{}{}\,(\P(r)\otimes\P(r))\otimes(A\otimes B)^{\otimes r} \cr
& \,\rto{\simeq}{}\,(\P(r)\otimes A^{\otimes r})\otimes(\P(r)\otimes B^{\otimes r}) \cr
& \,\rto{}{}\,A\otimes B \cr }$$
gives the tensor product $A\otimes B\in\Mod_\F$ the structure of a $\P$-algebra.
Equivalently,
let $\sum_i p_{(1)}^i\otimes p_{(2)}^i\in\P(r)\otimes\P(r)$ be the diagonal
of an operation $p\in\P(r)$.
Then, in $A\otimes B$, we have
$p(a_1\otimes b_1,\ldots,a_r\otimes b_r)
= \sum_i p_{(1)}^i(a_1,\ldots,a_r)\otimes p_{(2)}^i(b_1,\ldots,b_r)$.

\paragraph{\it The permutation operad}
The permutation groups $\Sigma_r$, $r\in\N$, form an operad in the category of sets.
This operad is called the permutation operad.
The composition product
$$\Sigma_r\times\Sigma_{s_1}\times\cdots\times\Sigma_{s_r}\,\rightarrow\,\Sigma_{s_1+\cdots+s_r}$$
is characterized by the relation
$$1_r(1_{s_1},\ldots,1_{s_r}) = 1_{s_1+\cdots+s_r}.$$
Thus, if $u\in\Sigma_r$ and $v_1\in\Sigma_{s_1},\ldots,v_r\in\Sigma_{s_r}$,
then, in $\Sigma_{s_1+\cdots+s_r}$,
we have the relation
$$\eqalign{ u(v_1,\ldots,v_r) &
= v_1\oplus\cdots\oplus v_r\cdot u_*(s_1,\ldots,s_r) \cr &
= u_*(s_1,\ldots,s_r)\cdot v_{u(1)}\oplus\cdots\oplus v_{u(r)}. \cr }$$
In fact,
the composite permutation $u(v_1,\ldots,v_r)\in\Sigma_{s_1+\cdots+s_r}$
can be obtained
by the following explicit substitution process.
In the sequence $(u(1),u(2),\ldots,u(r))$,
we replace the occurence of $k = 1,\ldots,r$
by the sequence $(v_k(1),v_k(2),\ldots,v_k(s_k))$
together with an appropriate shift.
To be precise,
we increase the elements $v_k(j)$ (where $j = 1,2,\ldots,s_k$)
by $s_1+s_2+\cdots+s_{k-1}$.
The process is the same for a partial composition product
$$u\circ_k v = u(1,\ldots,1,v,1,\ldots,1)\in\Sigma_{r+s-1}.$$
In this case, we have just to replace the occurence of $k$
by the sequence $(v(1)+k-1,v(2)+k-1,\ldots,v(s)+k-1)$
and to increase by $s$ the elements $u(i)$ such that $u(i)>k$.
For instance, if $u = (3, 2, 1)$ and $v = (1, 3, 2)$,
then we obtain $u\circ_2 v = (5, {\bf 2}, {\bf 4}, {\bf 3}, 1)$.

\paragraph{\it The associative and the commutative operads}\label{AssComOperads}
An algebra over the permutation operad is just an (associative) monoid.
There is also an operad in the category of sets whose algebras are commutative monoids
and whose components are reduced to base points $\O(r) = \pt$,
for all $r\in\N$.
We consider the associated operads in the category of $\F$-modules.

The associative operad $\A$ denotes the operad in the category of $\F$-modules
whose algebras are associative algebras.
The module $\A(r)$ is generated by the multilinear monomials in $r$ non-commutative variables.
Explicitly, we have:
$$\A(r) = \bigoplus\nolimits_{(i_1,\ldots,i_r)} \F\cdot x_{i_1}\cdots x_{i_r}$$
where $(i_1,\ldots,i_r)$ ranges over the permutations of $(1,\ldots,r)$.

The commutative operad $\C$ is the operad in the category of $\F$-modules
whose algebras are associative and commutative algebras.
The module $\C(r)$ is generated by the multilinear monomials in $r$ commutative variables.
Thus, we have
$$\C(r) = \F\cdot x_1\cdots x_r$$
and $\C(r)$ is the trivial representation of the symmetric group $\Sigma_r$.

There is an obvious operad morphism $\A\,\rto{}{}\,\C$.

\section{The operads}

\subsection{The Barratt-Eccles operad}

\paragraph{\it Summary}
In general,
an $E_\infty$-operad $\E$ is a $\Sigma_*$-projective resolution
of the commutative operad $\C$.
To be more explicit,
an operad is $\Sigma_*$-projective if it is projective as a $\Sigma_*$-module
(each component $\E(r)$ is a projective representation of $\Sigma_r$).
An operad is a resolution of the commutative operad $\C$,
if it is endowed with a quasi-isomorphism $\E\,\rto{\sim}{}\,\C$.

There is a canonical $E_\infty$-operad $\E$ given by the bar construction of the symmetric groups
and equipped with a diagonal $\E(r)\,\rightarrow\,\E(r)\otimes\E(r)$.
The aim of this section is to recall the definition
and to make explicit the structure of this operad.
In fact,
the operad $\E$ fits in a factorization of the morphism $\A\,\rto{}{}\,\C$
introduced in paragraph \ref{AssComOperads}.
To be more explicit, we have morphisms
$$\A\,\rto{}{}\,\E\,\rto{\sim}{}\,\C.$$
One observes precisely that the degree $0$ component of $\E(r)$
is equal to the associative operad $\A(r)$.
Accordingly,
we obtain an operad inclusion $\A\subset\E$.
Furthermore,
the maps $\E(r)_0 = \A(r)\,\rto{}{}\,\C(r)$
define an augmentation morphism $\E\,\rto{\sim}{}\,\C$.

In this article, the dg-operad $\E$ is called the {\it Barratt-Eccles operad},
because $\E$ is the normalized chain complex of a simplicial operad $\W$
introduced by M. Barratt and P. Eccles
for the study of infinite loop spaces ({\it cf}. [\cite{BE}]).

\paragraph{\it The dg-module $\E(r)$}
The dg-module $\E(r)$ is the normalized homogeneous bar complex
of the symmetric group $\Sigma_r$.
To be more explicit, the module $\E(r)_d$ is the $\F$-module
generated by the $d+1$-tuples $(w_0,\ldots,w_d)$,
where $w_i\in\Sigma_r$, for $i = 0,1,\ldots,d$.
If some consecutive permutations $w_i$ and $w_{i+1}$ are equal,
then, in $\E(r)_d$, we have $(w_0,\ldots,w_d) = 0$.
The differential of $\E(r)$ is given by the classical formula:
$$\delta(w_0,\ldots,w_r) = \sum_{i=0}^{d} (-1)^i (w_0,\ldots,\widehat{w_i},\ldots,w_d).$$
The permutations $\sigma\in\Sigma_r$ act diagonaly on $\E(r)_*$:
$$\sigma\cdot(w_0,\ldots,w_d) = (\sigma\cdot w_0,\ldots,\sigma\cdot w_d).$$

The quasi-isomorphism $\E(r)\,\rto{\sim}{}\,\C(r)$
is given by the classical augmentation $\E(r)\,\rto{}{}\,\F$
from the homogeneous bar construction to the trivial representation.
The degree $0$ component of $\E(r)$
is clearly the regular representation of the symmetric group $\Sigma_r$.
Therefore, we have a canonical morphism of dg-modules $\A(r)\,\rightarrow\,\E(r)$.

It is possible to specify an element $(w_0,\ldots,w_d)\in\E(r)$
by a table with $d+1$ rows indexed by $0,\ldots,d$.
The row $i$ consists of the sequence $(w_i(1),\ldots,w_i(r))$
and determines the permutation $w_i$.
For instance, the table:
$$\left|\matrix{ 1,2,3,4\hfill\cr 2,1,4,3\hfill\cr 1,4,2,3\hfill\cr }\right.$$
represents the element $(w_0,w_1,w_2)\in\E(4)_2$
such that $w_0 = (1,2,3,4)$, $w_1 = (2,1,4,3)$ and $w_2 = (1,4,2,3)$.

\paragraph{\it The composition product}
The composition product of the operad $\E$ is induced by the composition product
of permutations.
More precisely,
if $u = (u_0,\ldots,u_d)\in\E(r)_d$ and $v = (v_0,\ldots,v_e)\in\E(s)_e$,
then, in $\E(r+s-1)_{d+e}$,
we have the identity
$$u\circ_k v = \sum\nolimits_{(x_*,y_*)} \pm (u_{x_0}\circ_k v_{y_0},\ldots,u_{x_{d+e}}\circ_k v_{y_{d+e}}).$$
The sum ranges over the set of paths in an $d\times e$-diagram such as:
$$
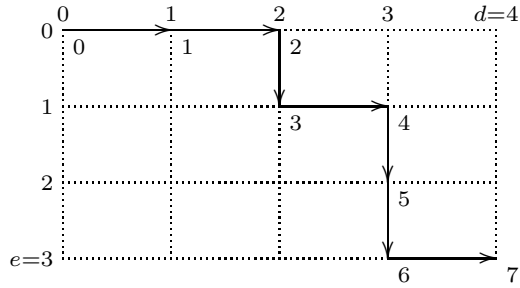
\xymatrix@M=0pt@C=12mm{\ar@{.}[d]_<{0}\ar@{->}[r]^<{0}_(0.15){0} &
\ar@{.}[d]\ar@{->}[r]^<{1}_(0.15){1} &
\ar@{->}[d]\ar@{.}[r]^<{2}_(0.15){2} &
\ar@{.}[d]\ar@{.}[r]^<{3}^>{d=4} &
\ar@{.}[d] \\
\ar@{.}[d]_<{1}\ar@{.}[r] & \ar@{.}[d]\ar@{.}[r] &
\ar@{.}[d]\ar@{->}[r]_(0.15){3} & \ar@{->}[d]\ar@{.}[r]_(0.15){4} & \ar@{.}[d] \\
\ar@{.}[d]_<{2}_>{e=3}\ar@{.}[r] & \ar@{.}[d]\ar@{.}[r] &
\ar@{.}[d]\ar@{.}[r] & \ar@{->}[d]\ar@{.}[r]_(0.15){5} & \ar@{.}[d] \\
\ar@{.}[r] & \ar@{.}[r] & \ar@{.}[r] & \ar@{->}[r]_(0.15){6}_(1.15){7} & \\ }$$
The indices involved in the sum are given by the coordinates of the vertices.
In the example above, we have:
$$(x_*,y_*)\,=\,(0,0),\,(1,0),\,(2,0),\,(2,1),\,(3,1),\,(3,2),\,(3,3),\,(4,3).$$
The sign associated to a term is determined by a permutation
of the horizontal and vertical segments of the path.
More precisely,
we consider the shuffle permutation
which takes the horizontal segments to the first places
and the vertical segments to the last places.
The sign is just the signature of this shuffle.
In the example, the segments
$$(0,1), (1,2), (2,3), (3,4), (4,5), (5,6), (6,7)$$
are permuted to
$$(0,1), (1,2), (3,4), (6,7), (2,3), (4,5), (5,6)$$
and, as a result, the sign is $-1$.

\paragraph{\it The diagonal}\label{AWDiag}
There is also a diagonal $\Delta: \E(r)\,\rightarrow\,\E(r)\otimes\E(r)$
which gives the cup-product
in the cohomology of the symmetric groups.
This diagonal is given by the classical formula:
$$\Delta(w_0,\ldots,w_n) = \sum_{d=0}^{n} (w_0,\ldots,w_d)\otimes(w_d,\ldots,w_n),$$
for any $(w_0,\ldots,w_n)\in\E(r)_d$.

\paragraph{\it On the classical Barratt-Eccles operad}
As mentioned in the introduction of the section,
the operad $\E$ is associated to a simplicial operad $\W$
(the classical Barratt-Eccles operad).
We recall the definition of $\W$.

If $X$ is a (discrete) set,
then $\W(X)$ denotes the classical contractible simplicial set associated to $X$,
which has the cartesian product
$$\W(X)_n = \underbrace{X\times\cdots\times X}_{n+1}$$
in dimension $n$.
The faces and degeneracies of $\W(X)$ are obtained by omitting and repeating components.
Explicitly, if $(x_0,\ldots,x_n)\in\W(X)_n$,
then we have
$$\eqalign{ d_i(x_0,\ldots,x_n) &
= (x_0,\ldots,\widehat{x_i},\ldots,x_n),
\qquad\hbox{for}\ i = 0,1,\ldots,n, \cr
\hbox{and}\qquad s_j(x_0,\ldots,x_n) &
= (x_0,\ldots,x_j,x_j,\ldots,x_n),
\qquad\hbox{for}\ j = 0,1,\ldots,n. \cr }$$
We observe that $X\mapsto\W(X)$ defines a (symmetric) monoidal functor from the category of sets
to the category of simplicial sets.
We have $\W(X\times Y) = \W(X)\times\W(Y)$.
As a consequence, if $\O$ is an operad in the category of sets,
then $\W(\O)$ is an operad in the category of simplicial sets.

The Barratt-Eccles operad is formed by the contractible simplicial sets
$\W(r) = \W(\Sigma_r)$, $r\in\N$, associated to the permutation operad.
To be explicit, the simplices $(w_0,\ldots,w_n)\in\W(\Sigma_r)_n$
support the diagonal action of the symmetric group.
We have:
$$\sigma\cdot(w_0,\ldots,w_n) = (\sigma\cdot w_0,\ldots,\sigma\cdot w_n),$$
for any $\sigma\in\Sigma_r$.
Furthermore, given $u = (u_0,\ldots,u_n)\in\W(\Sigma_r)_n$
and $v = (v_0,\ldots,v_n)\in\W(\Sigma_s)_n$,
the composite operation $u\circ_k v\in\W(\Sigma_{r+s-1})_n$ is given by the formula:
$$u\circ_k v = (u_0\circ_k v_0,\ldots,u_n\circ_k v_n).$$
The augmentation maps
$$\W(\Sigma_r)\,\rto{\sim}{}\,pt$$
form a morphism of simplicial operads,
since the base points define the operad of commutative monoids
({\it cf}. paragraph \ref{AssComOperads}).

We have clearly $\E(r) = N_*(\W(\Sigma_r))$.
In fact,
our definition of the composition product in the operad $\E$
follows from a classical representation of the Eilenberg-Zilber equivalence
$$N_d(\W(\Sigma_r))\otimes N_e(\W(\Sigma_s))\,\rto{\sim}{}\,N_{d+e}(\W(\Sigma_r)\times\W(\Sigma_s))$$
({\it cf}. [\cite{GZ}, Section II.5]).
The diagonal $\E(r)\,\rightarrow\,\E(r)\otimes\E(r)$ defined in paragraph \ref{AWDiag}
is clearly the Alexander-Whitney diagonal
of the normalized chain complex $N_*(\W(\Sigma_r))$.

\subsection{The surjection operad}

\smallskip
The {\it surjection operad} $\X$ is a quotient of the Barratt-Eccles operad $\E$.
The structure of the surjection operad $\X$ is made explicit in this section.
The quotient morphism is introduced in the next section.

\paragraph{\it The module of surjections}
The module $\X(r)_d$ is generated by the {\it non-degenerate} surjections
$u: \{1,\ldots,r+d\}\,\rightarrow\,\{1,\ldots,r\}$.
The degenerate surjections, whose definition is given below,
are equivalent to zero in $\X(r)_d$.
The non-surjective maps represent also the zero element in $X(r)_d$.
In the context of the surjection operad,
a map $u: \{1,\ldots,r+d\}\,\rightarrow\,\{1,\ldots,r\}$
is represented by a sequence $(u(1),\ldots,u(r+d))$.
By definition,
a surjection $u: \{1,\ldots,r+d\}\,\rightarrow\,\{1,\ldots,r\}$ is degenerate
if the associated sequence contains a repetition
(more explicitly, if we have $u(i) = u(i+1)$ for some $1\leq i<r+d$).
The module $\X(r)_d$ has a canonical $\Sigma_r$-module structure.
In fact, if $\sigma\in\Sigma_r$,
then the surjection $\sigma\cdot u\in\X(r)_d$ is represented
by the sequence $(\sigma(u(1)),\ldots,\sigma(u(r+d)))$.

We say that $v = (v(1),\ldots,v(s))$ is a {\it subsequence} of $u = (u(1),\ldots,u(r+d))$
if $v(1) = u(i_1), v(2) = u(i_2), \ldots, v(s) = u(i_s)$,
for some $1\leq i_1<i_2<\cdots<i_s\leq r+d$.

\paragraph{\it The table arrangement of a surjection}
We define the {\it table arrangement} of a surjection $u\in\X(r)_d$.
The table arrangement is obtained by splitting $u$ into certain subsequences
$$u_i = (\underbrace{u(r_0+\cdots+r_{i-1}+1)}_{u_i(1)},
\underbrace{u(r_0+\cdots+r_{i-1}+2)}_{u_i(2)},\ldots,
\underbrace{u(r_0+\cdots+r_{i-1}+r_i)}_{u_i(r_i)})$$
where $r_0,r_1,\ldots,r_d\geq 1$ and $r_0+r_1+\cdots+r_d = r+d$.
In general,
the values of $u$ are just arranged on a table with $d+1$ rows indexed by $0,\ldots,d$.
The subsequence $u_i$ consists of the terms of the row $i$:
$$u(1),\ldots,u(r+d)\,=\,\left|\matrix{ u_0(1),\hfill & \cdots & u_0(r_0-1),\hfill & u_0(r_0),\hfill \cr
u_1(1),\hfill & \cdots & u_1(r_1-1),\hfill & u_1(r_1),\hfill \cr
\quad\vdots\hfill & & \quad\vdots\hfill & \quad\vdots\hfill \cr
u_{d-1}(1),\hfill & \cdots & u_{d-1}(r_{d-1}-1),\hfill & u_{d-1}(r_{d-1}),\hfill \cr
u_d(1),\hfill & \cdots & u_d(r_d).\hfill \cr }\right.$$
Accordingly,
the subsequence $u_i$ is referred to as the line $i$ of the table arrangement of $u$.
The elements
$$u_0(r_0),u_1(r_1),\ldots,u_{d-1}(r_{d-1})$$
which are at the end of a row (except for the last row)
are called the {\it caesuras} of the surjection.

The caesuras are the elements of the sequence $(u(1),\ldots,u(r+d))$
which do not represent the last occurence of a value $k = 1,\ldots,r$.
Accordingly,
the other elements
$$u_0(1),\ldots,u_0(r_0-1),\ldots,u_{d-1}(1),\ldots,u_{d-1}(r_{d-1}-1),u_d(1),\ldots,u_d(r_d),$$
represent all the last occurence of a value in the sequence $(u(1),\ldots,u(r+d))$.
(In particular, these elements form a permutation of $1,\ldots,r$.)
Observe that the number of caesuras equals the degree $d$ of the surjection.

As an example, the sequence $u = (1,3,2,1,4,2,1)$ has the following table arrangement:
$$1, 3, 2, 1, 4, 2, 1\,=\,\left|\matrix{ 1,\hfill\cr 3, 2,\hfill\cr 1,\hfill\cr 4, 2, 1.\hfill\cr }\right.$$

\paragraph{\it The complex of surjections}\label{SurjectionDiff}
In this paragraph,
we equip the module of surjections with a differential
$\delta: \X(r)_*\,\rightarrow\,\X(r)_{*-1}$.
The identity $\delta^2 = 0$ follows from results
which we prove in the next section.

By definition,
if a surjection $u\in\X(r)_d$ is specified by a sequence $(u(1),\ldots,u(r+d))$,
then the differential $\delta(u)\in\X(r)_{d-1}$ equals the sum:
$$\delta(u(1),\ldots,u(r+d)) = \sum_{i=1}^{r+d} \pm (u(1),\ldots,\widehat{u(i)},\ldots,u(r+d)).$$
The signs are determined by the table arrangement.
To be precise,
we mark each term of the surjection $u(i)$ with a sign
that we insert in the formula of the differential.
The caesuras of the surjection are marked with alternate signs.
If a value occurs only once,
then we observe that the associated sequence $(u(1),\ldots,\widehat{u(i)},\ldots,u(r+d))$
does not represent a surjection and vanishes in $\X(r)_{d-1}$.
Otherwise, the last occurence of a value is marked with the sign opposite
to the sign associated to the previous occurence of the same value
(which is a caesura).

As an example,
in case of the surjection $u = (1,3,2,1,4,2,1)\in\X(4)_3$,
we obtain the signs:
$$\left|\matrix{ 1^+,\hfill\cr 3, 2^-,\hfill\cr 1^+,\hfill\cr 4, 2^+, 1^-.\hfill\cr }\right.$$
Hence, we have $\delta(1,3,2,1,4,2,1) = (3,2,1,4,2,1)
- (1,3,1,4,2,1) + (1,3,2,4,2,1) + (1,3,2,1,4,1) - (1,3,2,1,4,2)$.

\paragraph{\it The operadic composition of surjections}\label{CompositeSurjectn}
We equip the module of surjections with a partial composition product
$\X(r)_d\otimes\X(s)_e\,\rightarrow\,\X(r+s-1)_{d+e}$.
If $u = (u(1),\ldots,u(r+d))\in\X(r)_d$ and $v = (v(1),\ldots,v(s+e))\in\X(s)_e$,
then the product $u\circ_k v\in\X(r+s-1)_{d+e}$
is obtained by substituting the occurences of $k$ in $(u(1),\ldots,u(r+d))$
by elements of $(v(1),\ldots,v(s+e))$.

Precisely, assume that $k$ has $n$ occurences in $(u(1),\ldots,u(r+d))$
which are the terms $u(i_1),\ldots,u(i_n)$.
In this case, we split the sequence $(v(1),\ldots,v(s+e))$ into $n$ com\-po\-nents:
$$(v(j_0),\ldots,v(j_1))\quad(v(j_1),\ldots,v(j_2))\quad\cdots\quad(v(j_{n-1}),\ldots,v(j_n))$$
where $1 = j_0\leq j_1\leq j_2\leq\cdots\leq j_{r-1}\leq j_n = s+e$.
Then, in $(u(1),\ldots,u(r+d))$,
we replace the term $u(i_m)$ by the sequence $(v(j_{m-1}),\ldots,v(j_m))$.
In addition, we increase the terms $v(j)$ by $k-1$.
The terms $u(i)$ such that $u(i)<k$ are fixed.
The terms $u(i)$ such that $u(i)>k$ are increased by $s-1$.
We call {\it composite surjection} any sequence which arises from this process.
The composition product $u\circ_k v$ is represented by the sum of all composite surjections
together with a sign which we specify in the next paragraph.

As an example, we have:
$$(1,2,1,3)\circ_1(1,2,1) = \pm(1,3,1,2,1,4)\pm(1,2,3,2,1,4)\pm(1,2,1,3,1,4).$$
In fact, there are 2 occurences of $k = 1$ in $(1,2,1,3)$.
Therefore, we cut the sequence $(1,2,1)$ in 2 pieces.
There are 3 possibilities $(1) (1,2,1)$, $(1,2) (2,1)$ and $(1,2,1) (1)$.
Hence, the substitution gives the composite surjections
$({\bf 1},3,{\bf 1},{\bf 2},{\bf 1},4)$,  $({\bf 1},{\bf 2},3,{\bf 2},{\bf 1},4)$
and $({\bf 1},{\bf 2},{\bf 1},3,{\bf 1},4)$.

\paragraph{\it The sign associated to a composite surjection}\label{SignCompositeSurjectn}
The sign associated to a composite surjection is determined by sequence permutations
which are involved in our substitution process.

In more details, the subsitution process works as follows.
The sequence $(u(1),\ldots,u(r+d))$ is cut at the occurences of $k$.
Thus, with the notation above, we obtain the components:
$$(u(i_0),\ldots,u(i_1))
\quad(u(i_1),\ldots,u(i_2))
\quad\cdots\quad(u(i_{n-1}),\ldots,u(i_n))
\quad(u(i_n),\ldots,u(i_{n+1})),$$
where, by convention, $i_0 = 1$ and $i_{n+1} = r+d$.
The sequence $(v(1),\ldots,v(s+e))$ is split into
$$(v(j_0),\ldots,v(j_1))\quad(v(j_1),\ldots,v(j_2))\quad\cdots\quad(v(j_{n-1}),\ldots,v(j_n)).$$
We insert the sequence $(v(j_{m-1}),\ldots,v(j_m))$ at the position of $u(i_m)$.
Hence,
performing compositions from right to left,
we have to permute this sequence with the components
$$(u(i_m),\ldots,u(i_{m+1}))\quad\cdots\quad(u(i_n),\ldots,u(i_{n+1})).$$
The sequences are equipped with a degree (defined in the next paragraph).
Therefore,
according to the sign-rule,
these permutations determine a sign
which affects the result of the substitution process.

By definition, the sequence $(u(i_{m-1}),\ldots,u(i_m))$ has degree $d'$
if it intersects $d'+1$ lines
in the table arrangement of the surjection $(u(1),\ldots,u(r+d))$.
Accordingly,
if $d_m$ denotes the degree of a sequence $(u(i_{m-1}),\ldots,u(i_m))$,
then $(u(i_{m-1}),\ldots,u(i_m))$ intersects the lines
$$d_1+\cdots+d_{m-1},d_1+\cdots+d_{m-1}+1,\ldots,d_1+\cdots+d_{m-1}+d_m.$$
The degree of the components of $(v(1),\ldots,v(s+e))$ are given by the same rule.

\paragraph{\it Example}
As an example,
let us determine the signs which occur in the composition product
$$(1,2,1,3)\circ_1(1,2,1) = \pm(1,3,1,2,1,4)\pm(1,2,3,2,1,4)\pm(1,2,1,3,1,4)$$
considered above.
The sequence $(1,2,1,3)$ has the following table arrangement:
$$1,2,1,3\,=\,\left|\matrix{ 1,\hfill\cr 2, 1, 3.\hfill\cr }\right.$$
This sequence is split into $(1)(1,2,1)(1,3)$.
The component $(1)$ which is contained in the first line has degree 0.
The component $(1,2,1)$ intersects the lines 0 and 1 and has degree 1.
The component $(1,3)$ contained in the last line has degree 0.
The degrees of the components of $(1,2,1)$ are determined similarly.
Now, we perform the following permutations and concatenations
(the subscripts denote the degree of the components)
$$\eqalign{ & (1)_{0}(1,3,1)_{1}(1,4)_{0}({\bf 1})_{0}({\bf 1},{\bf 2},{\bf 1})_{1}
\mapsto {\bf +}(1)_{0}({\bf 1})_{0}(1,3,1)_{1}({\bf 1},{\bf 2},{\bf 1})_{1}(1,4)_{0}
\mapsto {\bf +}({\bf 1},3,{\bf 1},{\bf 2},{\bf 1},4), \cr
& (1)_{0}(1,3,1)_{1}(1,4)_{0}({\bf 1},{\bf 2})_{1}({\bf 2},{\bf 1})_{0}
\mapsto {\bf -}(1)_{0}({\bf 1},{\bf 2})_{1}(1,3,1)_{1}({\bf 2},{\bf 1})_{0}(1,4)_{0}
\mapsto {\bf -}({\bf 1},{\bf 2},3,{\bf 2},{\bf 1},4), \cr
& (1)_{0}(1,3,1)_{1}(1,4)_{0}({\bf 1},{\bf 2},{\bf 1})_{1}({\bf 1})_{0}
\mapsto {\bf -}(1)_{0}({\bf 1},{\bf 2},{\bf 1})_{1}(1,3,1)_{1}({\bf 1})_{0}(1,4)_{0}
\mapsto {\bf -}({\bf 1},{\bf 2},{\bf 1},3,{\bf 1},4), \cr }$$
from which we deduce the sign associated to each term
of the product $(1,2,1,3)\circ_1(1,2,1)$.
We obtain finally
$(1,2,1,3)\circ_1(1,2,1) = (1,3,1,2,1,4) - (1,2,3,2,1,4) - (1,2,1,3,1,4)$.

\smallskip
We claim that our definitions make sense:

\paragraph{\sc Proposition}\label{SurjectionOpd}
{\it The constructions above provide the graded module of surjections
with the structure of a dg-operad.}

\smallskip
Explicitly, we prove that the composition product of surjections
is compatible with the differential
and satistifies the unit, associativity and equivariance properties
of an operad composition product.
In fact, we deduce this proposition from results
which we prove in the next sections.
The surjection operad is also a $\Sigma_*$-projective resolution
of the commutative operad ({\it cf}. [\cite{MCSSeq}]).
But, we do not use this result in the article.

\subsection{The table reduction morphism}

\paragraph{\it The table reduction morphism}
We define a morphism $\TR: \E\,\rightarrow\,\X$ (the {\it table reduction morphism}).
Let $w = (w_0,\ldots,w_d)\in\E(r)_d$.
Given indices $r_0,\ldots,r_d\geq 1$ such that $r_0+\cdots+r_d = r+d$,
we form a surjection $w' = (w'(1),\ldots,w'(r+d))\in\X(r)_d$
from the values of the permutations $w_0,\ldots,w_d$.
We say that this surjection $w'\in\X(r)_d$ arises from the simplex $w\in\E(r)_d$
by a table reduction process.
The image of $w = (w_0,\ldots,w_d)$ in $\X(r)$ is the sum of these elements
$$\TR(w_0,\ldots,w_d) = \sum\nolimits_{(r_0,\ldots,r_d)} (w'(1),\ldots,w'(r+d))$$
for all choices of $r_0,\ldots,r_d\geq 1$.
Let us mention that there are only positive signs in the sum.

Fix $r_0,\ldots,r_d\geq 1$ as above.
The associated surjection has a table arrangement
$$w'(1),\ldots,w'(r+d)\,=\,\left|\matrix{ w'_0(1),\hfill & \cdots & w'_0(r_0-1),\hfill & w'_0(r_0),\hfill \cr
w'_1(1),\hfill & \cdots & w'_1(r_1-1),\hfill & w'_1(r_1),\hfill \cr
\quad\vdots\hfill & & \quad\vdots\hfill & \quad\vdots\hfill \cr
w'_{d-1}(1),\hfill & \cdots & w'_{d-1}(r_{d-1}-1),\hfill & w'_{d-1}(r_{d-1}),\hfill \cr
w'_d(1),\hfill & \cdots & w'_d(r_d),\hfill \cr }\right.$$
whose rows have length $r_0,r_1,\ldots,r_d$
and are given by subsequences of the permutations $w_0,\ldots,w_d$.

Precisely, the table can be obtained as follows.
The elements $(w'_0(1),\ldots,w'_0(r_0))$ are the first $r_0$ terms
of the permutation $(w_0(1),\ldots,w_0(r))$.
The elements $(w'_i(1),\ldots,w'_i(r_i))$ are the first $r_i$ terms
of the permutation $(w_i(1),\ldots,w_i(r))$
which do not occur in
$$\matrix{ w'_0(1),\hfill & \cdots & w'_0(r_0-1),\hfill\cr
w'_1(1),\hfill & \cdots & w'_1(r_1-1),\hfill\cr
\quad\vdots\hfill & & \quad\vdots\hfill\cr
w'_{i-1}(1),\hfill & \cdots & w'_{i-1}(r_{i-1}-1).\hfill\cr }$$
(These elements are supposed to represent the final occurence of a value in the surjection
and cannot be repeated below in the table.)

For instance, we have:
$$\TR\left|\matrix{ 1, 2, 3, 4\hfill\cr 1, 4, 3, 2\hfill\cr 1, 2, 4, 3\hfill\cr }\right.
\,=\,\left|\matrix{ 1, 2,\hfill\cr 4,\hfill\cr 2, 4, 3\hfill\cr }\right.
\,+\,\left|\matrix{ 1, 2,\hfill\cr 4, 3,\hfill\cr 2, 3.\hfill\cr }\right.$$
(We remove the degenerate terms, which are equivalent to $0$, from the result.)

\paragraph{\sc Theorem}\label{TRPpty}
{\it The table reduction morphism $\TR: \E\,\rightarrow\,\X$
defined in the paragraph above
is a surjective morphism of differential graded operads.}

\smallskip
We prove proposition \ref{SurjectionOpd} and theorem \ref{TRPpty} together.
Here is the plan of our demonstration.
First, we have the following results:

\paragraph{\sc Assertion:}\label{TRNormalized}
{\it The table reduction morphism $\TR: \E\,\rightarrow\,\X$ is well defined.}

\paragraph{\sc Lemma}\label{TREpi}
{\it The table reduction morphism $\TR: \E\,\rightarrow\,\X$ is surjective.}

\smallskip
Then:

\paragraph{\sc Assertion:}\label{TRDiff}
{\it The table reduction morphism $\TR: \E\,\rightarrow\,\X$
maps the differential of the Barratt-Eccles operad $\E(r)$ to the differential
of the surjection operad $\X(r)$ (specified in paragraph \ref{SurjectionDiff}).}

\smallskip
This implies that the differential of $\X(r)$ verifies the identity $\delta^2 = 0$.

\smallskip
The next assertion is obvious from the definition of the table reduction morphism:

\paragraph{\sc Assertion:}
{\it The table reduction morphism $\TR: \E\,\rightarrow\,\X$ is equivariant.}

\smallskip
The next assertion implies that the composition product of the surjection operad
satisfies the properties of an operad composition product
and is compatible with the differential.

\paragraph{\sc Assertion:}\label{TRComp}
{\it The table reduction morphism $\TR: \E\,\rightarrow\,\X$
maps the composition product in the Barratt-Eccles operad to the composition product
in the surjection operad (specified in paragraph \ref{CompositeSurjectn}).}

\smallskip
These statements prove that the surjection operad $\X$ is a quotient operad
of the Barratt-Eccles operad $\E$.

\smallskip
Assertions \ref{TRNormalized}, \ref{TREpi} and \ref{TRDiff}
are proved in section \ref{TRPPtyProofs}.
Section \ref{TRCompProofs} is devoted to the proof
of the last assertion \ref{TRComp}.

\subsection{The operad morphism properties}\label{TRPPtyProofs}

\paragraph{\it Proof of assertion \ref{TRNormalized}:}
If $w\in\E(r)_d$ is a degenerate simplex in the Barratt-Eccles operad,
then the associated surjections $w'\in\X(r)_d$ are also degenerate.
Therefore, the element $\TR(w)$ vanishes in the normalized surjection operad.
This proves that the table reduction morphism $\TR: \E\,\rightarrow\,\X$ is well defined.
To be more explicit,
we assume $w = s_j(w_0,\ldots,w_d) = (w_0,\ldots,w_j,w_j,\ldots,w_d)$.
Consider a surjection $w'$ which arises from the table reduction
of $(w_0,\ldots,w_j,w_j,\ldots,w_d)$.
In the associated table arrangement,
the line $j$ consists of the terms
$$w'_j(1),\ldots,w'_j(s_j-1),w'_j(s_j)$$
of the permutation $w_j$.
The values $w'_j(1),\ldots,w'_j(s_j-1)$ are omitted in the next row
(by definition of the table reduction process).
Therefore, this row (which consists also of terms of the permutation $w_j$)
starts with the term $w'_j(s_j)$.
We conclude that the last element of the line $j$
coincides with the first element of the line $j+1$.
The assertion follows.

\paragraph{\it Proof of lemma \ref{TREpi}:}\label{TRSectn}
Let $u\in\X(r)_d$ be a surjection.
We specify a simplex $w = (w_0,w_1,\ldots,w_d)\in\E(r)_d$ such that $\TR(w) = u$.
The permutations $(w_0,w_1,\ldots,w_d)$ are subsequences of $(u(1),\ldots,u(r+d))$ 
and are determined inductively from the table arrangement of the surjection.

To be explicit, the table arrangement of $u$ has the form:
$$u(1),\ldots,u(r+d)\,=\,\left|\matrix{ s(1),\hfill & \cdots & s(i_1-1),\hfill & x(1),\hfill \cr
s(i_1),\hfill & \cdots & s(i_2-1),\hfill & x(2),\hfill \cr
\quad\vdots\hfill & & \quad\vdots\hfill & \quad\vdots\hfill \cr
s(i_{d-1}),\hfill & \cdots & s(i_d-1),\hfill & x(d),\hfill \cr
s(i_d),\hfill & \cdots & s(r),\hfill \cr }\right.$$
where $(s(1),\ldots,s(r))$ is a permutation.
We fix $w_d = (s(1),\ldots,s(r))$.
In order to obtain the permutation $w_i$ from $w_{i+1}$,
we just move the occurence of $x = x(i+1)$ in the given subsequence
$$(w_{i+1}(1),\ldots,w_{i+1}(r))\subset(u(1),\ldots,u(r+d))$$
to the position of the caesura $x(i+1)$.
Explicitly, we have:
$$\matrix{ w_d\hfill & = s(1),\ldots,s(r),\hfill\cr
w_{d-1}\hfill & = s(1),\ldots,s(i_d-1),x(d),s(i_d),\ldots,\widehat{x(d)},\ldots,s(r),\hfill\cr
w_{d-2}\hfill & = s(1),\ldots,s(i_{d-1}-1),x(d-1),s(i_{d-1}),\ldots,\widehat{x(d-1)},\ldots\qquad\hfill\cr
& \hfill\ldots,s(i_d-1),x(d),s(i_d),\ldots,\widehat{x(d)},\ldots,s(r)\cr }$$
and so on.
The lemma follows from the following claim:

\paragraph{\sc Claim:} {\it If $(w_0,\ldots,w_d)$ are defined as above,
then we have $\TR(w_0,\ldots,w_d) = (u(1),\ldots,u(r+d))$.}

\smallskip
The permutations $w_i$ and $w_{i+1}$ coincide up to the position of $x(i)$.
To be more explicit:
$$\matrix{ w_i\hfill & = w_i(1),\ldots,w_i(p-1),x(i+1),w_i(p+1),\ldots,w_i(q),w_i(q+1),\ldots,w_i(r),\hfill\cr
w_{i+1}\hfill & = w_i(1),\ldots,w_i(p-1),w_i(p+1),\ldots,w_i(q),x(i+1),w_i(q+1),\ldots,w_i(r).\hfill\cr }$$
Now, consider a term
$$\left|\matrix{ w'_0(1),\hfill & \cdots & w'_0(r_0-1),\hfill & w'_0(r_0),\hfill \cr
w'_1(1),\hfill & \cdots & w'_1(r_1-1),\hfill & w'_1(r_1),\hfill \cr
\quad\vdots\hfill & & \quad\vdots\hfill & \quad\vdots\hfill \cr
w'_{d-1}(1),\hfill & \cdots & w'_{d-1}(r_{d-1}-1),\hfill & w'_{d-1}(r_{d-1}),\hfill \cr
w'_d(1),\hfill & \cdots & w'_d(r_d)\hfill \cr }\right.$$
in $\TR(w_0,\ldots,w_d)$.
The caesura in the permutation $w_i$ is necessarily put at $x(i+1)$.
Otherwise, in the table, there would be a repetition between the lines $i$ and $i+1$.
Precisely, one observes that the last term of the line $i$ would coincide with the first term of the line $i+1$.
This property holds for any permutation $w_i$.
The claim follows.

\paragraph{\it Proof of assertion \ref{TRDiff}:}
Fix $w = (w_0,\ldots,w_d)\in\E(r)_d$.
Let $w'\in\X(r)_d$ be a surjection associated to $w$ by the table reduction process.
If we assume $r_i = 1$, then we have:
$$w'(1),\ldots,w'(r+d)\,=\,\left|\matrix{ w'_0(1),\ldots,w'_0(r_0),\hfill\cr
\quad\vdots\hfill\cr
w'_{i-1}(1),\ldots,w'_{i-1}(r_{i-1}),\hfill\cr
w'_i(1),\hfill\cr
w'_{i+1}(1),\ldots,w'_{i+1}(r_{i+1}),\hfill\cr
\quad\vdots\hfill\cr
w'_d(1),\ldots,w'_d(r_d).\hfill\cr }\right.$$
Clearly, in the differential $\delta(w')\in\X(r)_{d-1}$,
the term, which removes $w'_i(1)$ from the sequence,
is equivalent to a surjection associated to $(w_0,\ldots,\widehat{w_i},\ldots,w_d)\in\E(r)_{d-1}$
by the table reduction process.
Furthermore, in both cases, the sign in the differential is $(-1)^i$.
We observe next that the other terms in $\delta(\TR(w_0,\ldots,w_d))$ cancel each other.
Therefore, the differential $\delta(\TR(w_0,\ldots,w_d))\in\X(r)_{d-1}$
reduces to $\TR(\delta(w)) = \sum_i\pm\TR(w_0,\ldots,\widehat{w_i},\ldots,w_d)$.

To be precise, we observe that the omission of a final occurence in a surjection $w'\in\X(r)_d$
cancels the omission of a caesura in another surjection  $w''\in\X(r)_d$ associated to $w$.
Explicitly, let $w'\in\X(r)_d$ be the surjection
specified by the indices $r_0,\ldots,r_d$
in the table reduction of $w\in\E(r)_d$.
Assume that we remove the element $w'_i(p)$ from the sequence $w'(1),\ldots,w'(r+d)$.
Assume that this element $w'_i(p)$ is the final occurence of the value $k$ in $w'$.
As in the definition of the differential,
we consider the previous occurence of $k$ in $w'$,
which is at a caesura $w'_j(r_j)$ in the table arrangement of the surjection.
Let $w''\in\X(r)_d$ be the surjection specified
by the parameters $r_0,\ldots,r_j+1,\ldots,r_i-1,\ldots,r_d$
in the table reduction of $w\in\E(r)_d$.
This surjection $w''$ differs from $w'$ by the lines $i$ and $j$
of the associated table arrangement:
$$w''\,=\,\left|\matrix{ w'_0(1),\ldots,w'_0(r_0),\hfill\cr
\quad\vdots\hfill\cr
w'_j(1),\ldots,w'_j(r_j-1),k,w'_j(r_j+1),\hfill\cr
\quad\vdots\hfill\cr
w'_i(1),\ldots,w'_i(p-1),\widehat{k},w'_i(p+1),\ldots,w'_i(r_i),\hfill\cr
\quad\vdots\hfill\cr
w'_d(1),\ldots,w'_d(r_d).\hfill\cr }\right.$$
We observe that the omission of $w'_j(r_j+1)$ in $w''$ yields the same result
as the omission of $w'_i(p)$ in $w'$.
Furthermore, the signs associated to these differentials are opposite
(this assertion is an immediate consequence of the definition).
The conclusion follows.

\subsection{On composition products}\label{TRCompProofs}

\smallskip
The purpose of this section is to prove assertion \ref{TRComp}
and to complete the proofs of proposition \ref{SurjectionOpd}
and theorem \ref{TRPpty}.
We compare the composition product of the surjection operad
with the composition product of the Barratt-Eccles operad.
To be precise, we prove the following lemma:

\paragraph{\sc Lemma}
{\it Let $u'\in\X(r)_d$ (respectively, $v'\in\X(r)_e$) be a surjection
which comes from the table reduction of a simplex $u\in\E(r)_d$ (respectively, $v\in\E(s)_e$)
in the Barratt-Eccles operad.

1) For each non-degenerate surjection $w'\in\X(r+s-1)_{d+e}$
in the expansion of $u'\circ_k v'\in\X(r+s-1)_{d+e}$,
there is one and only one $d\times e$-path $(x_*,y_*)$
such that $w'$ is obtained by table reduction
from $(u_{x_*}\circ_k v_{y_*})\in\E(r+s-1)_{d+e}$.

2) Conversely, any non-degenerate surjection obtained by table reduction from $u\circ_k v$
arises from a composite surjection $w'$ as in 1).

3) The signs associated to the surjection $w'$ in the expansion of $\TR(u)\circ_k\TR(v)$
and in the expansion of $\TR(u\circ_k v)$ coincide.}

\smallskip
By definition of the table reduction morphism,
the surjections associated to the simplex $(u_{x_*}\circ_k v_{y_*})\in\E(r+s-1)_{d+e}$
are all distinct
because the associated table arrangements have distinct dimensions.
Therefore, the statements 1), 2), 3) permit to identify the non-degenerate terms
in the expansion of $\TR(u)\circ_k\TR(v)$
with the non-degenerate terms in the expansion of $\TR(u\circ_k v)$.
Consequently, we have $\TR(u)\circ_k\TR(v) = \TR(u\circ_k v)$,
as claimed by assertion \ref{TRComp}.

\smallskip
We specify the path $(x_*,y_*)$ associated to the surjection $w'$
in paragraph \ref{PathCompSurjection}.
We prove that $w'$ arises from $(u_{x_*}\circ_k v_{y_*})$ by table reduction
in claim \ref{TRCompProcess}.
The assertion 3) in the lemma is immediate from the construction of the path
({\it cf}. observation \ref{TRCompSign}).
The assertion 2) follows from claim \ref{TRCompSurjective}.
Finally, in part 1), the uniqueness property is a consequence of a stronger result
(claim \ref{TRCompInjective}).

\paragraph{\it On the table arrangement of a composite surjection}
We need insights in the substitution process of surjections.
We fix $u'\in\X(r)_d$ and $v'\in\X(s)_e$ as in the lemma above.
By convention, in the sequel, the lines of the table arrangement of $u'$ (respectively, $v'$)
are indexed by $x = 0,\ldots,d$ (respectively, $y = 0,\ldots,e$)
and have length $r_0,\ldots,r_d$ (respectively, $s_0,\ldots,s_e$).

As in paragraphs \ref{CompositeSurjectn} and \ref{SignCompositeSurjectn},
we let $u'(i_1),\ldots,u'(i_n)$ denote the occurences of $k$
in the sequence $u'(1),\ldots,u'(r+d)$.
By convention, we have also $i_0=1$ and $i_{n+1}=r+d$.
The surjection $u'$ is split into
$$(u'(i_0),\seg,u'(i_1))\quad(u'(i_1),\seg,u'(i_2))\quad\cdots\quad(u'(i_n),\seg,u'(i_{n+1})).$$
Furthermore, we fix a splitting
$$(v'(j_0),\seg,v'(j_1))\quad(v'(j_1),\seg,v'(j_2))\quad\cdots\quad(v'(j_{n-1}),\seg,v'(j_n))$$
of the surjection $v'$.
We insert these components in the splitting of $u'$
$$\displaylines{ (u'(i_0),\seg,u'(i_1))\quad(v'(j_0),\seg,v'(j_1))
\quad(u'(i_1),\seg,u'(i_2))\quad\cdots\hfill\cr
\hfill\cdots\quad(v'(j_{n-1}),\seg,v'(j_n))\quad(u'(i_n),\seg,u'(i_{n+1})) \cr }$$
and we let
$$\displaylines{ w' = (u'(i_0),\seg,u'(i_1-1),v'(j_0),\seg,v'(j_1),u'(i_1+1),\ldots\hfill\cr
\hfill\ldots,u'(i_{n-1}-1),v'(j_{n-1}),\seg,v'(j_n),u'(i_n+1),\seg,u'(i_{n+1})).\cr }$$
For simplicity,
we shall omit the substitution-shift in this section.
For instance, in the formula above, we omit to increase the elements $v'(j)$ by $k-1$
and to increase the elements $u'(i)>k$ by $s-1$.

\smallskip
Each caesura of the composite surjection
comes either from a caesura of $u'$ or from a caesura of $v'$.
More precisely:

\paragraph{\sc Observation:}\label{CaesurasCompSurjection}
{\it The caesuras of the composite surjection $w'$ are classified as follows.
Each caesura $u'_x(r_x)$ of the surjection $u'$ such that $u'_x(r_x)\not=k$
gives a caesura in the composite surjection.
Each caesura $u'_x(r_x)$ such that $u'_x(r_x)=k$
gives a caesura in the composite surjection
at the last element $v'(j_m)$ of the sequence $v'(j_{m-1}),\seg,v'(j_m)$
which is inserted at $u'_x(r_x)$.
The other elements of this sequence which are caesuras of $v'$
give caesuras in the composite surjection.}

\smallskip
Our next purpose is to relate the table arrangement of the composite surjection
to the composition procedure in the Barratt-Eccles operad.

\paragraph{\it The path determined by the insertion of a surjection}\label{PathCompSurjection}
We define the path $(x_*,y_*)$ associated to a composite surjection.
We let $d_m$ and $e_m$ denote the respective degrees of $(u'(i_{m-1}),\seg,u'(i_m))$
and $(v'(j_{m-1}),\seg,v'(j_m))$.
The definition implies that the sequence $(u'(i_{m-1}),\seg,u'(i_m))$
intersects the lines
$$d_1+\cdots+d_{m-1},d_1+\cdots+d_{m-1}+1,\ldots,d_1+\cdots+d_{m-1}+d_m$$
in the table arrangement of $u'$.
Similarly, the sequence $(v'(i_{m-1}),\seg,v'(i_m))$ intersects the lines
$$e_1+\cdots+e_{m-1},e_1+\cdots+e_{m-1}+1,\ldots,e_1+\cdots+e_{m-1}+e_m$$
in the table arrangement of $v'$.

By definition, the path $(x_*,y_*)$ associated to the composite surjection
is given by the concatenation of horizontal components
of length $d_m$, $m = 1,\ldots,n+1$,
and of vertical components of length $e_m$, $m = 1,\ldots,n$,
as in the insertion process.
The indices $(x_c,y_c)$ are given by the following formulas.
$$\eqalign{ & \hbox{If}\ c = d_1+e_1+d_2+e_2+\cdots+d_{m-1}+e_{m-1}+x',\ \hbox{where}\ 0\leq x'\leq d_m, \cr
& \hbox{then we have}\ \left\{\matrix{\ x_c = d_1+d_2+\cdots+d_{m-1}+x',\hfill\cr
\ y_c = e_1+e_2+\cdots+e_{m-1}.\hfill\cr }\right. \cr
& \hbox{If}\ c = d_1+e_1+d_2+e_2+\cdots+d_{m-1}+e_{m-1}+d_m+y',\ \hbox{where}\ 0\leq y'\leq e_m, \cr
& \hbox{then we have}\ \left\{\matrix{\ x_c = d_1+d_2+\cdots+d_{m-1}+d_m,\hfill\cr
\ y_c = e_1+e_2+\cdots+e_{m-1}+y'.\hfill\cr }\right. \cr }$$
This definition is motivated by the table arrangement of the composite surjection.

\paragraph{\sc Observation:}
{\it In the component $(u'(i_{m-1}),\seg,u'(i_m))$,
the elements $u'(i)\not=k$,
which come from the line $x = d_1+\cdots+d_{m-1}+x'$ of the table arrangement of $u'$,
are inserted in the line $c = d_1+e_1+\cdots+d_{m-1}+e_{m-1}+x'$
of the table arrangement of the composite surjection $w'$.

In the component $(v'(j_{m-1}),\seg,v'(j_m))$,
the elements $v'(j)$
which come from the line $y = e_1+\cdots+e_{m-1}+y'$ of the table arrangement of $v'$,
are inserted in the line $c = d_1+e_1+\cdots+d_{m-1}+e_{m-1}+d_m+y'$
of the table arrangement of the composite surjection $w'$.}

\smallskip
This observation follows from observation \ref{CaesurasCompSurjection}.

\smallskip
It is easy to compare the sign produced by the composition
of the surjections $u'$ and $v'$
to the sign associated to the path $(x_*,y_*)$.
In the composition of $u'$ and $v'$,
the sign is produced by the insertion of the components $(v'(j_{m-1}),\ldots,v'(j_m))$
of the surjection $v'$
in the surjection $u'$.
The sign associated to the path $(x_*,y_*)$
is determined by a shuffle of the horizontal and vertical components
of the path.
It is immediate from the definitions that we consider the same shuffles.
Therefore:

\paragraph{\sc Observation:}\label{TRCompSign}
{\it The sign produced by the insertion of a surjection agrees
with the sign determined by the path $(x_*,y_*)$
associated to the composite surjection.}

\smallskip
Now, we have:

\paragraph{\sc Claim:}\label{TRCompProcess}
{\it The composite surjection $w'$ can be obtained by table reduction of the simplex
$(u_{x_*}\circ_k v_{y_*})\in\E(r+s-1)_*$,
where $(x_*,y_*)$ is the path associated to $w'$.}

\proof
We let $(u\circ_k v)'$ be the surjection in the expansion of $\TR(u_{x_*}\circ_k v_{y_*})$
whose table arrangement has the same dimensions as the table arrangement of $w'$.
We prove by induction that $(u\circ_k v)'$ agrees with $w'$.
Given $0\leq x\leq d$, we assume that the table arrangement of $(u\circ_k v)'$
agrees with the table arrangement of $w'$ on the lines $c$
such that $x_c<x$.
We prove that the lines $c = c',c'+1,\ldots,c''$ such that $x_c = x$ agree also.
We have:
$$(x_{c'},y_{c'}) = (x,y'),\ (x_{c'+1},y_{c'+1}) = (x,y'+1),\ \ldots,\ (x_{c''},y_{c''}) = (x,y'').$$

In the table arrangement of the composite surjection $w'$,
the lines $c',c'+1,\ldots,c''$ either do not contain or contain an occurence of the surjection $v'$.
In the first case, we have $c' = c''$ and the line considered reads:
$$\left|\matrix{ u'_x(1),\ldots,u'_x(r_x).\hfill\cr }\right.$$
In the second case, the lines $c',c'+1,\ldots,c''$ of the table arrangement
of the composite surjection have the following form:
$$\left|\matrix{ u'_x(1),\ldots,u'_x(i'-1),v'_{y'}(j'),\ldots,v'_{y'}(s_{y'}),\hfill\cr
v'_{y'+1}(1),\ldots,v'_{y'+1}(s_{y'+1}),\hfill\cr
\qquad\vdots\hfill\cr
v'_{y''-1}(1),\ldots,v'_{y''-1}(s_{y''-1}),\hfill\cr
v'_{y''}(1),\ldots,v'_{y''}(s_{y''}),u'_x(i'+1),\ldots,u'_x(r_x),\hfill\cr }\right.$$
where $1\leq i'\leq r_x$.
It is also possible to have $c' = c''$.
Then, the diagram above collapses to one line:
$$\left|\matrix{ u'_x(1),\ldots,u'_x(i'-1),
v'_y(j'),\ldots,v'_{y'}(s_{y'}),
u'_x(i'+1),\ldots,u'_x(r_x),\hfill\cr }\right.$$

In the table arrangement of the surjection $(u\circ_k v)'$,
the lines $c = c',c'+1,\ldots,c''$ arise
from the table reduction of the permutations $u_x\circ_k v_y$,
where $y = y', y'+1, \ldots, y''$.
Explicitly, the permutation $u_x\circ_k v_y$ is represented by the sequence:
$$u_x(1),\ldots,u_x(p-1),v_y(1),\ldots,v_y(s),u_x(p+1),\ldots,u_x(r),$$
where $u_x(p) = k$.
We have just to omit the elements which already occur inside a line above
in the table arrangement of $(u\circ_k v)'$.
By the induction hypothesis,
these elements occupy the same position
in the table arrangement of $w'$.
Therefore, the elements which we have to removed from the permutations
$u_x\circ_k v_y$, $y = y', y'+1, \ldots, y''$
are the same as the elements which are omitted
in the lines $c = c',c'+1,\ldots,c''$
of the composite surjection $w'$.

Let us follow the definition of the table reduction morphism for $(u\circ_k v)'$.
Let us assume that we have not to omit all elements $v_{y'}(1),\ldots,v_{y'}(s)$
in the construction of the line $c'$ of $(u\circ_k v)'$.
As observed, the elements $u_x(q)$ which we have to remove from $u_x(1),\ldots,u_x(p-1)$
are also omitted in the sequence $u'_x(1),\ldots,u'_x(r_x)$.
It follows that the first terms inserted on the line $c'$ of the surjection $(u\circ_k v)'$
and of the composite surjection $w'$ agree.
If the line is not completed,
then the elements
which we insert from $v_{y'}(1),\ldots,v_{y'}(s)$
and from $u_x(p+1),\ldots,u_x(r)$
are also the same
for the same reason.
If $c'<c''$, then on the lines $c = c'+1,\ldots,c''$
of the table arrangement of $(u\circ_k v)'$,
we omit $u_x(1),\ldots,u_x(p-1)$
because these elements occur either inside the line $c'$
or inside a line above in the table.
Therefore, the next elements that we have to insert
come from $v_y(1),\ldots,v_y(s)$, $y = y'+1,\ldots,y''$.
These elements correspond again to the next elements $v'_y(1),\ldots,v'_y(s_{y'})$
of the composite surjection $w'$.

If we have to omit all elements $v_{y'}(1),\ldots,v_{y'}(s)$,
then, for the table arrangement of the composite surjection,
we obtain:
$$\left|\matrix{ u'_x(1),\ldots,u'_x(r_x).\hfill\cr }\right.$$
Furthermore, the elements $u'_x(1),\ldots,u'_x(r_x)$ belong necessarily
to the last component $(u'(i_n),\ldots,u'(i_{n+1}))$
of the surjection $u'$.
Therefore, the value $k$ has to be omitted from $u_x(1),\ldots,u_x(p-1),k,u_x(p+1),\ldots,u_x(r)$
in the construction of $u'_x(1),\ldots,u'_x(r_x)$.
The omission of the block $v_y(1),\ldots,v_y(s)$ in the permutation $u_x\circ_k v_y$ gives the same result.
The conclusion follows in this case.

The proof of the claim is complete.\endproof

\paragraph{\sc Claim:}\label{TRCompSurjective}
{\it Fix a surjection $(u\circ_k v)'$ in the expansion of $\TR(u\circ_k v)$.
Fix $x = 0,\ldots,d$.
We assume $x_c = x$, for $c = c',c'+1,\ldots,c''$.
Consider the sequence formed by the lines $c = c',c'+1,\ldots,c''$ of the surjection $(u\circ_k v)'$:
$$\left|\matrix{ (u_x\circ_k v_{y'})'(1),\ldots,(u_x\circ_k v_{y'})'(t_{c'}),\hfill\cr
(u_x\circ_k v_{y'+1})'(1),\ldots,(u_x\circ_k v_{y'+1})'(t_{c'+1}),\hfill\cr
\qquad\vdots\hfill\cr
(u_x\circ_k v_{y''})'(1),\ldots,(u_x\circ_k v_{y''})'(t_{c''}).\hfill\cr }\right.$$
This sequence is equivalent to a concatenation:
$$(u'_x(1),\ldots,u'_x(i'-1),k)\quad(v'(j'),\ldots,v'(j''))\quad(k,u'_x(i'+1),\ldots,u'_x(r_x)),$$
otherwise the surjection $(u\circ_k v)'$ is degenerate.}

\proof
The property is immediate if $c'=c''$.
We assume $c'<c''$.
In the table, all lines except the last one have to end with a term $v'_y(j)$
and all lines except the first one have to start with a term $v'_y(1)$,
otherwise the table contains a repetition.

In fact, if the line $y$ ends with $u'_x(i')$, then all values $u_x(i)$,
which come before $u'_x(i')$ in the sequence $u_x(1),\ldots,\widehat{k},\ldots,u_x(r)$,
are omitted in the next lines of the surjection.
If $k$ comes before $u'_x(i')$,
then all values $v_{y}(1),\ldots,v_{y}(s)$ occur already
as a final element in the table.
Therefore, the elements $v_{y}(1),\ldots,v_{y}(s)$ are also omitted in the next lines.
As a conclusion, the line $y+1$ has to start with the element $u'_x(i')$
and a repetition occurs in the table.

Now, if we assume that the line $y-1$ ends with $v'_{y-1}(j)$,
then the elements $u_x(i)$ which come before $k$ occur already as a final element in the table.
Therefore, the line $y$ starts with an element $v'_{y}(1)$.

The claim follows.\endproof

\paragraph{\sc Claim:}\label{TRCompInjective}
{\it Let $w'\in\X(r+s-1)$ be any surjection
which occurs in the expansion of $u'\circ_k v'$
for some $u'\in\X(r)$ and $v'\in\X(s)$.
This composite surjection $w'$ is non-degenerate
if and only if the surjections $u'$ and $v'$ are non-degenerate.
In this case, the surjections $u'\in\X(r)$ and $v'\in\X(s)$
and the splitting of $v'$ which gives rise to $w'$
are uniquely determined by $w'\in\X(r+s-1)$.}

\proof
The first assertion is immediate from the insertion process of surjections.
The degeneracies $u'(i) = u'(i+1)$ such that $u'(i)\not=k$ are not removed by the insertion process.
So do the degeneracies $v'(j) = v'(j+1)$.
If the surjection $u'$ has a degeneracy $u'(i) = u'(i+1)$, such that $u'(i) = k$,
then, according to our notation, we have $i = i_m$ and $i+1 = i_{m+1}$ ,
for some $m = 0,1,\ldots,n$.
In the composite surjection,
this degeneracy is equivalent to the repetition of the element $v'(j_m)$.

The sequence $u'(1),\ldots,u'(r)$ is recovered by the following operations.
The elements such that $k\leq w'(i)\leq k+t-1$ are replaced by $k$.
The elements such that $k+t\leq w'(i)$ are decreased by $t-1$.
The consecutive occurences of $k$ are collapsed to one term
and this achieves the process.
(If the surjection $u'$ is non-degenerate,
then the repetition of $k$ follows from the insertion process.)
Similarly, for the sequence $v'(1),\ldots,v'(s)$,
we withdraw the elements such that $w'(i)\leq k-1$ or $k+t\leq w'(i)$
and we decrease by $k-1$ the elements such that $k\leq w'(i)\leq k+t-1$.
Consecutive occurences of a same value are also collapsed to one term.

This achieves the proof of the claim.\endproof

\smallskip
The proof of assertion \ref{TRComp} is now complete.

\subsection{On little cubes filtrations}

\smallskip
The Barratt-Eccles operad has a filtration
$$F_1\W(r)\subseteq F_2\W(r)\subseteq\cdots\subseteq F_n\W(r)\subseteq\cdots\subseteq\W(r)$$
by simplicial suboperads $F_n\W(r)$
whose topological realizations $|F_n\W|$ are homotopy equivalent
to the classical operads of little $n$-cubes $F_n\D$
({\it cf}. [M. Boardmann, R. Vogt, \cite{BV}] and [P. May, \cite{MLoop}]).
Precisely,
the topological operads $|F_n\W|$ and $F_n\D$ are connected by operad morphisms
$$|F_n\W|\,\lto{\sim}{}\,\cdot\,\rto{\sim}{}\,F_n\D$$
which are homotopy equivalences of topological spaces.

The filtration of the Barratt-Eccles operad is introduced by J. Smith
in the article [\cite{JfS}]
in order to generalize the Milnor $FK$-construction
and to provide a simplicial model for iterated loop spaces.
The work of T. Kashiwabara ({\it cf}. [\cite{K}])
proves that the simplicial set $F_n\W(r)$
has the same homology as the little $n$-cubes space $F_n\D(r)$.
The equivalence as a topological operad
follows from the work of the first author ({\it cf}. [\cite{Bg}]).
The idea is to relate the filtration to a particular cellular structure.

We have an induced filtration on the associated dg-operad
$$F_1\E(r)\subseteq F_2\E(r)\subseteq\cdots\subseteq F_n\E(r)\subseteq\cdots\subseteq\E(r)$$
such that $F_n\E(r)$ is equivalent to the dg-operad
formed by the chain complexes of the little $n$-cubes spaces.
J. McClure and J. Smith have introduced a similar filtration
$$F_1\X(r)\subseteq F_2\X(r)\subseteq\cdots\subseteq F_n\X(r)\subseteq\cdots\subseteq\X(r)$$
for the surjection operad ({\it cf}. [\cite{MCSSeq}]).
Furthermore,
these authors prove that the dg-operad $F_n\X$ is equivalent to the little cubes operad.
This implies that we have quasi-isomorphisms:
$$F_n\E\,\lto{\sim}{}\,\cdot\,\rto{\sim}{}\,F_n\X.$$
Our purpose is to prove the following lemma:

\paragraph{\sc Lemma}\label{CubeFiltratn}
{\it The table reduction morphism $\TR: \E\,\rightarrow\,\X$ preserves filtrations.
Moreover, the induced morphisms $F_n\TR: F_n\E\,\rightarrow\,F_n\X$
are quasi-isomorphisms of differential graded operads.}

\paragraph{\it The cellular structures}
We recall the definition of the cellular structures of the operads $\P = \E$ and $\P = \X$.
A cell is associated to a pair $(\mu,\sigma)$
such that $\sigma\in\Sigma_r$ is a permutation
and $\mu$ is a collection of non-negative integers $\mu_{i j}\in\N$
indexed by pairs $i<j$,
where $i,j\in\{1,\ldots,r\}$.
The permutation $\sigma\in\Sigma_r$
determines also a collection of permutations
$\sigma_{i j}\in\Sigma_{\{i,j\}}$.
Explicitly,
the element $\sigma_{i j}$ is the permutation of $(i,j)$
formed by the occurences of $i$ and $j$
in the sequence $(\sigma(1),\ldots,\sigma(r))$
which represents the permutation $\sigma$.
Equivalently,
$\sigma_{i j}$ is the identity permutation $(i,j)$
if $\sigma^{-1}(i)<\sigma^{-1}(j)$
and $\sigma_{i j}$ is the transposition $(j,i)$
if $\sigma^{-1}(i)>\sigma^{-1}(j)$.

The cellular structure is specified by sub-$\Sigma_r$-modules
$F_{(\mu,\sigma)}\P(r)\subseteq\P(r)$.
The filtration is related to the cellular structure by the relation:
$$F_n\P(r) = \sum_{\mu_{i j}<n} F_{(\mu,\sigma)}\P(r).$$
The sum ranges over all pairs $(\mu,\sigma)$
such that $\mu_{i j}<n$ for all $i<j$.

The cellular structure is in some sense compatible
with the operad composition product ({\it cf}. [\cite{Bg}]).
Furthermore,
the cells $F_{(\mu,\sigma)}\P(r)$ are acyclic
in both cases $\P = \E$ ({\it cf}. [\cite{Bg}])
and $\P = \X$ ({\it cf}. [\cite{MCSSeq}]).
In the next paragraphs,
we just recall the definition of the cells $F_{(\mu,\sigma)}\P(r)\subseteq\P(r)$.

\paragraph{\it The cellular structure of the Barratt-Eccles operad}
Let $w = (w_0,\ldots,w_d)\in\E(r)_d$.
Given a pair $i<j$, we consider the sequence
$w_{i j} = ((w_0)_{i j},\ldots,(w_d)_{i j})$
formed by the permutations of $i$ and $j$
associated to $(w_0,\ldots,w_d)$.
The simplex $w = (w_0,\ldots,w_d)$ belongs to $F_{(\mu,\sigma)}\E(r)_d\subset\E(r)_d$
if, for all pairs $i<j$,
either the sequence $w_{i j}$ has no more than $\mu_{i j}-1$ variations,
or the sequence $w_{i j}$ has exactly $\mu_{i j}$ variations
and we have $(w_d)_{i j} = \sigma_{i j}$.
Accordingly,
the simplex $w = (w_0,\ldots,w_d)$ belongs to $F_n\E(r)_d\subset\E(r)_d$
if and only if, for all pairs $i<j$, the sequence $w_{i j}$
has no more than $n-1$ variations.
For instance,
we have $((1,2),(2,1),(1,2))\in F_3\W(2)_2$.

\paragraph{\it The cellular structure of the surjection operad}
Fix a surjection $u\in\X(r)_d$.
Given a pair $i<j$, we let $u_{i j}$ denote the subsequence of $(u(1),\ldots,u(r+d))$
formed by the occurences of $i$ and $j$ in the surjection.
For instance, we have $(2,1,3,4,3,1)_{1\,2} = (2,1,1)$,
$(2,1,3,4,3,1)_{1\,3} = (1,3,3,1)$ and $(2,1,3,4,3,1)_{2\,3} = (2,3,3)$.
The surjection $u\in\X(r)_n$ belongs to $F_{(\mu,\sigma)}\X(r)_d\subset\X(r)_d$
if, for all pairs $i<j$,
either the sequence $u_{i j}$ has no more than $\mu_{i j}$ variations,
or the sequence $u_{i j}$ has exactly $\mu_{i j}+1$ variations
and the permutation formed by the final occurences of $i$ and $j$ in $u_{i j}$
agrees with $\sigma_{i j}$.
Accordingly,
the surjection $u\in\X(r)_n$ belongs to $F_n\X(r)_d\subset\X(r)_d$
if and only if, for all pairs $i<j$,
the sequence $u_{i j}$ has no more than $n$ variations.
For instance, we have $(1,2,1,2)\in F_3\X(2)_2$.
Similarly, for the surjection above, we have $(2,1,3,4,3,2)\in F_2\X(4)_2$.

\smallskip
Lemma \ref{CubeFiltratn} follows from the following result:

\paragraph{\sc Lemma}\label{CellFiltratn}
{\it The table reduction morphism $\TR: \E\,\rightarrow\,\X$ preserves the cellular structures.
Explicitly, if $w\in F_{(\mu,\sigma)}\E$, then we have $\TR(w)\in F_{(\mu,\sigma)}\X$.}

\smallskip
The result above implies that the induced morphism $F_n\TR: F_n\E\,\rightarrow\,F_n\X$
is a weak equivalence, because we have the commutative diagram:
$$\xymatrix{ *{\displaystyle\hocolim_{(\mu,\sigma)} F_{(\mu,\sigma)}\E}
\ar[rr]^{F_{(\mu,\sigma)}\TR}\ar[d]^{\sim} & &
*{\displaystyle\hocolim_{(\mu,\sigma)} F_{(\mu,\sigma)}\X}
\ar[d]^{\sim} \\
F_n\E\ar[rr]^{F_n\TR} & & F_n\X \\ }$$
({\it cf}. [\cite{Bg}], [\cite{MCSSeq}]).
The upper horizontal arrow is a weak equivalence
because the cells $F_{(\mu,\sigma)}\E$ and $F_{(\mu,\sigma)}\X$
are acyclic.

In the note [\cite{BF}], we define a chain morphism $\TC: \X\,\rightarrow\,\E$
quasi-inverse to $\TR: \E\,\rightarrow\,\X$
and such that $\TC(F_n\X)\subset F_n\E$.
This construction gives another proof of lemma \ref{CubeFiltratn}.

The proof of the lemma above is postponed to the end of this section.

\paragraph{\it The relationship to operations on the Hochschild cochains}
According to J. McClure and J. Smith,
the operad $F_2\X$ is generated by the operations $(1,2)\in\X(2)$
and $(1,2,1,3,1,\ldots,1,r,1)\in\X(r)$.
Furthermore, if $A$ is an associative algebra,
then the normalized Hochschild cochain complex $C^*(A,A)$
is equipped with the structure of an algebra over $F_2\X$
({\it cf}. [\cite{MCSPrism}], [\cite{MCSSeq}]).
The generator $(1,2)\in\X(2)$ is mapped to the cup-product operation
$$x_1\smile x_2: C^*(A,A)\otimes C^*(A,A)\,\rightarrow\,C^*(A,A).$$
The generators $(1,2,1,3,1,\ldots,1,r,1)\in\X(r)$ are mapped to the brace-operations
$$x_1\{x_2,\ldots,x_r\}: C^*(A,A)^{\otimes r}\,\rightarrow\,C^*(A,A)$$
introduced by E. Getzler ({\it cf}. [\cite{Ge}])
and by T. Kadeishvili ({\it cf}. [\cite{Kd}]).

J. McClure and J. Smith conclude that the surjection operad
is an instance of a dg-operad
which is equivalent to the little squares operad
and which operates on the normalized Hochschild cochain complexes.
This result answers a conjecture of P. Deligne.
Lemma \ref{CubeFiltratn} implies the same result about our simplicial model of the little squares operad
(regardless of the homotopy type of the operad $F_2\X$).
Moreover, claim \ref{TRSectn} allows to define elements $w\in F_2\E$
which give rise to the cup-product of Hochschild cochains
and to the Getzler-Kadeishvili brackets.
Finally, we obtain the following result:

\paragraph{\sc Theorem}
{\it The normalized Hochschild cochain complex of an associative algebra $C^*(A,A)$
is equipped with the structure of an algebra over the operad $F_2\E$.
The cup-product of Hoch\-schild cochains
$$x_1\cup x_2: C^*(A,A)\otimes C^*(A,A)\,\rightarrow\,C^*(A,A)$$
is the operation associated to the element $w = (\id)\in\E(2)_0$.
The Getzler-Kadeishvili bracket
$$x_1\{x_2,\ldots,x_r\}: C^*(A,A)^{\otimes r}\,\rightarrow\,C^*(A,A)$$ 
is the operation associated to the element $w = (w_0,\ldots,w_{r-1})\in\E(r)_{r-1}$
such that
$$\eqalign{ & w_0 = (1,2,3,4,\ldots,r-1,r),\quad w_1 = (2,1,3,4,\ldots,r-1,r), \cr
& w_2 = (2,3,1,4,\ldots,r-1,r),\quad\ldots\quad w_{r-1} = (2,3,4,\ldots,r-1,r,1). \cr }$$}

\paragraph{\it Proof of lemma \ref{CellFiltratn}}
Let $w = (w_0,w_1,\ldots,w_d)\in\E(r)_d$.
Let $w'$ be a surjection associated to $w$ by the table reduction process.
Fix a pair $i<j$.
If $w_{i j}$ has no more than $n$ variations,
then we claim that $w'_{i j}$ has no more than $n+1$ variations.

We assume that the final occurence of $i$
comes before the final occurence of $i$ in the surjection $w'$
(the arguments are the same in the inverse case).
We assume that the final occurence of $i$ in $w'$
belongs to the line $m$ of the table arrangement of $w'$.
Consider the subsequence of $w'_{i j}$
which is formed by all the occurences of $i$ in the surjection $w'$
and by the occurences of $y$ which are on the lines $0,1,\ldots,m-1$
of the table arrangement of $w'$.
Equivalently, the sequence $w'_{i j}$ is truncated at the final occurence of $i$.
In $w'_{i j}$, there is at least an occurence of $j$
which comes after the final occurence of $i$.
Therefore, the sequence $w'_{i j}$ has one more variation than the truncated subsequence.

If $i$ occurs in the line $k$ of the table arrangement,
then $i$ comes necessarily before $j$ in the sequence $w_k(1),\ldots,w_k(r_k)$,
because we assume that the final occurence of $j$ lies below in the table.
Hence, in this case, we have $(w_k)_{i j} = (i, j)$.
Similarly, if $j$ occurs in a line $w'_k(1),\ldots,w'_k(r_k)$, such that $k<m$,
then $j$ comes necessarily before $i$ in the sequence $w_k(1),\ldots,w_k(r)$.
Hence, in this case, we have $(w_k)_{i j} = (j, i)$.
These observations imply that the truncated subsequence introduced above has no more variations
than $(w_0)_{i j},(w_1)_{i j},\ldots,(w_m)_{i j}$.
The conclusion follows.

If $w'_{i j}$ has exactly one more variation than $w_{i j}$,
then the permutation formed by the final occurences of $i$ and $j$ in $w'$
agree with $(w_d)_{i j}$.
In fact, under this assumption,
the sequence $w_{i j} = ((w_0)_{i j},(w_1)_{i j},\ldots,(w_d)_{i j})$
has as much variation as the truncated sequence
$((w_0)_{i j},(w_1)_{i j},\ldots,(w_m)_{i j})$.
Therefore, we have $(w_m)_{i j} = (w_{m+1})_{i j} = \cdots = (w_d)_{i j}$.
The conclusion follows from the observations above
(we have $(w_m)_{i j} = (i,j)$ in the case considered).

The lemma is proved.

\paragraph{\it Remark}
It is also straighforward to prove that the induced morphism
$F_{(\mu,\sigma)}\TR: F_{(\mu,\sigma)}\E\,\rightarrow\,F_{(\mu,\sigma)}\X$
is surjective.
To be precise, given $u\in\X(r)_d$, we consider the simplex $w\in\E(r)_d$
introduced in paragraph \ref{TRSectn} and such that $\TR(w) = u$.
>From the definition, it follows easily that $w\in F_{(\mu,\sigma)}\E$
as long as $u\in F_{(\mu,\sigma)}\X$.

We can improve the construction of the simplex $w\in\E(r)_d$.
More precisely, we have a morphism of dg-modules $\TC: \X\,\rightarrow\,\E$
which preserves the cellular structures (but not the operad composition products)
and which is right inverse to the table reduction morphism $\TR: \E\,\rightarrow\,\X$
({\it cf}. [\cite{BF}]).

\section{The interval cut operations on chains}

\parnb=-1\paragraph{\it Conventions for simplicial sets}\label{ConvtnSimplicial}
We refer to D. Curtis ({\it cf}. [\cite{C}]),
P. Gabriel and M. Zisman ({\it cf}. [\cite{GZ}]),
P. Goerss and J. Jardine ({\it cf}. [\cite{GoJ}]),
P. May ({\it cf}. [\cite{MSimp}])
for our conventions on simplicial sets.
We adopt the classical notation $[n]$, $n\in\N$,
for the objects of the simplicial category $\Delta$.
A morphism $u: [m]\,\rightarrow\,[n]$ is equivalent to a sequence
$0\leq u(0)\leq u(1)\leq\cdots\leq u(m)\leq n$.
If $X$ is a simplicial set and $x\in X_n$ is an $n$-dimensional simplex of $X$,
then the expression $x(u(0),\ldots,u(m))\in X_m$ denotes the image of $x\in X_n$
under the map $u^*: X_n\,\rightarrow\,X_m$ associated to $u: [m]\,\rightarrow\,[n]$.
For instance, we have the identity $x = x(0,\ldots,n)$.
According to this convention,
we have also:
$$\eqalign{ d_i(x) & = x(0,\ldots,i-1,i+1,\ldots,n),
\quad\hbox{for}\ i = 0,1,\ldots,n,\cr
\hbox{and}\quad s_j(x) & = x(0,\ldots,j,j,\ldots,n),
\quad\hbox{for}\ j = 0,1,\ldots,n.\cr }$$

We consider the universal simplicial set $\Delta^n$ such that $(\Delta^n)_m = \Hom_{\Delta}([m],[n])$.
In regard to the conventions above,
the simplex of $\Delta^m$ associated to $u: [m]\,\rightarrow\,[n]$
is denoted by $\Delta(u(0),\ldots,u(m))\in(\Delta^n)_m$.
Given $x\in X_n$,
we let $\tilde{x}: \Delta^n\,\rightarrow\,X$
denote the unique morphism such that $\tilde{x}(\Delta(0,\ldots,n)) = x$.
We have in fact $\tilde{x}(\Delta(u(0),\ldots,u(m))) = x(u(0),\ldots,u(m))$.
Furthermore,
if $u: [m]\,\rightarrow\,[n]$,
then $u_*: \Delta^m\,\rightarrow\,\Delta^n$ is the unique morphism
which maps the simplex $\Delta(0,\ldots,m)\in(\Delta^m)_m$
to $\Delta(u(0),\ldots,u(m))\in(\Delta^n)_m$.

\subsection{The Eilenberg-Zilber operad}

\smallskip
The purpose of this section is to prove the following theorem:

\paragraph{\sc Theorem}\label{BarrattEcclesChainOperations}
{\it Let $X$ be a simplicial set. We have natural evaluation coproducts
$\E(r)\otimes N_*(X)\,\rightarrow\,N_*(X)^{\otimes r}$
which make the normalized chain complex $N_*(X)$
into a coalgebra over the Barratt-Eccles operad $\E$.
Accordingly,
the dual cochain complex $N^*(X)$ is equipped with the structure of an $\E$-algebra.

We let $\Sigma_2 = \{\id,\tau\}$ where $\id = (1,2)$ and $\tau = (2,1)$.
The cooperation $\mu_0^*: N_*(X)\,\rightarrow\,N_*(X)^{\otimes 2}$
associated to the element $\mu_0 = (\id)\in\E(2)_0$
coincides with the classical Alexan\-der-Whitney diagonal.
The dual operation $\mu_0: N^*(X)^{\otimes 2}\,\rightarrow\,N^*(X)$
is nothing but the classical cup-product of cochains.
The higher operation $\mu_d: N^*(X)^{\otimes 2}\,\rightarrow\,N^*(X)$
associated to the sequence $\mu_d = (\id,\tau,\id,\tau,\ldots)\in\E(2)_d$
is a representative of the cup-$d$-product.}

\paragraph{\it On coalgebras over an operad}
An algebra over an operad $\P$ is defined as a dg-module $A$
together with an evaluation product
$\P(r)\otimes A^{\otimes r}\,\rightarrow\,A$.
Dually,
a coalgebra over an operad $\P$ is a dg-module $C$
together with an equivariant evaluation coproduct
$\P(r)\otimes C\,\rightarrow\,C^{\otimes r}$
which is coassociative with respect to the operad product
and unital with respect to the operad unit.
This definition implies clearly that the dual dg-module of a $\P$-coalgebra
is equipped with the structure of a $\P$-algebra.

In the context of coalgebras,
an element $p\in\P(r)$ determines a cooperation $p^*: C\,\rightarrow\,C^{\otimes r}$.
In fact,
the dg-modules $\Hom^c_C(r) = \Hom_\F(C,C^{\otimes r})$ form a dg-operad
and the structure of a $\P$-coalgebra is equivalent
to a morphism of dg-operads $\P\,\rightarrow\,\Hom^c_C$.
To be more explicit, let $V$ be a dg-module.
The endomorphism operad of $V$ is the dg-operad such that $\Hom_V(r) = \Hom_\F(V^{\otimes r},V)$
(see paragraph \ref{EndOp}).
Dually, the {\it coendomorphism operad} associated to $V$ is the dg-operad
such that $\Hom^c_V(r) = \Hom_\F(V,V^{\otimes r})$.
Thus, if $g_1: V\,\rightarrow\,V^{\otimes s_1}$,
\dots, $g_r: V\,\rightarrow\,V^{\otimes s_r}$
and $f: V\,\rightarrow\,V^{\otimes r}$,
then $f(g_1,\ldots,g_r): V\,\rightarrow\,V^{\otimes s_1+\cdots+s_r}$
is the composite morphism
$f(g_1,\ldots,g_r) = \pm g_1\otimes\cdots\otimes g_r\cdot f$.
The sign follows from the permutation of $f$ with $g_1,\ldots,g_r$.
The dg-module $V$ is equipped with the structure of a coalgebra
over the coendomorphism operad $\Hom^c_V$
and the canonical morphism
$\Hom^c_V(r)\otimes V\,\rightarrow\,V^{\otimes r}$
is a universal evaluation coproduct.

\paragraph{\it The Eilenberg-Zilber operad}
The Eilenberg-Zilber operad $\Z$ is the universal dg-operad
together with a natural evaluation coproduct
$$\Z(r)\otimes N_*(X)\,\rightarrow\,N_*(X)^{\otimes r}$$
defined for $X\in\sSet$.
We refer to V. Hinich and V. Schechtman ({\it cf}. [\cite{HS}])
and V. Smirnov ({\it cf}. [\cite{Sm}]).
The idea of the construction goes back to the work of A. Dold
on Steenrod operations ({\it cf}. [\cite{D}]).
The dg-module $\Z(r)$ is formed by the morphisms
$$f_X: N_*(X)\,\rightarrow\,N_*(X)^{\otimes r}$$
in the product of the coendomorphism operads
$\prod_X\Hom_\F(N_*(X),N_*(X)^{\otimes r})$
which define a natural transformation in $X\in\sSet$.
Therefore, for a fixed $X\in\sSet$,
we have a canonical operad morphism
$\Z(r)\,\rightarrow\,\Hom_\F(N_*(X),N_*(X)^{\otimes r})$
which is equivalent to an evaluation coproduct as above.

Classically, a natural transformation $f_X: C_*(X)\,\rightarrow\,N_*(X)^{\otimes r}$
is determined by the image of the simplices
$\Delta(0,\ldots,n)\in C_n(\Delta^n)$,
where $n\in\N$.
Furthermore, a natural morphism $f_X: C_*(X)\,\rightarrow\,N_*(X)^{\otimes r}$
induces a natural transformation on the normalized complex of $X$
if the tensors $f_X(\Delta(0,\ldots,n))\in N_*(\Delta^n)^{\otimes r}$, $n\in\N$,
are cancelled by the degeneracies $s^j: \Delta^n\,\rightarrow\,\Delta^{n+1}$.
To conclude, we have the identity:
$$\Z(r) = \prod_n\biggl\{\bigcap_{j=0}^n
\ker(s^j:\,N_*(\Delta^n)^{\otimes r}\,\rightarrow\,N_*(\Delta^{n+1})^{\otimes r})
\biggr\}.$$
To be more precise,
a tensor of degree $d$ in $N_*(\Delta^n)^{\otimes r}$
is associated to an operation
of degree $d-n$ in the Eilenberg-Zilber operad.

\subsection{The construction of interval cut operations}

\smallskip
We construct an operad morphism $\AW: \X\,\rightarrow\,\Z$.
We have then a sequence of operad morphisms
$$\E\,\rto{\TR}{}\,\X\,\rto{\AW}{}\,\Z,$$
which, according to the discussion above,
provide the chain complexes $N_*(X)$, $X\in\sSet$,
with the structure of a coalgebra over $\E$,
as claimed by theorem \ref{BarrattEcclesChainOperations}.
The cooperation $\AW(u): N_*(X)\,\rightarrow\,N_*(X)^{\otimes r}$,
which is associated to a surjection $u\in\X(r)_d$,
generalizes the Alexander-Whitney diagonal.
This justifies the notation $\AW: \X\,\rightarrow\,\Z$
of our morphism.
The operation $\AW(u)\in\Z(r)_d$ associated to a surjection $u\in\X(r)_d$
is also called an {\it interval cut operation}
because the simplices $\Delta(C_{(1)}),\ldots,\Delta(C_{(r)})$
which occur in the expansion of $\AW(u)(\Delta(0,\ldots,n))\in N_*(\Delta^n)^{\otimes r}$
arise by splitting the interval $\{0,1,\ldots,n\}$.
The details of the construction are given in the next paragraphs.

Our generalization of the Alexander-Whitney diagonal goes back to the original construction
of the reduced square operations by N. Steenrod ({\it cf}. [\cite{S}]).
A morphism $\AW: \X\,\rightarrow\,\Z$ similar to ours
is defined by D. Benson ({\it cf}. [\cite{Bn}, Section 4.5])
and J. McClure and J. Smith ({\it cf}. [\cite{MCSSeq}]).
Our construction differs simply by sign conventions.

\paragraph{\it On the interval cut associated to a surjection}\label{IntervalCut}
We define the simplices $\Delta(C_{(k)})\in N_*(\Delta^n)$, $k = 1,\ldots,r$,
which are associated to a surjection $u\in\X(r)_d$
for a fixed splitting of the interval $\{0,1,\ldots,n\}$.
In fact, we have $\Delta(C_{(k)}) = c_{(k)}^*\Delta(0,\ldots,n)$
for certain morphisms $c_{(k)}: [m_k]\,\rightarrow\,[n]$, $k = 1,\ldots,r$,
in the simplicial category.
Thus, the letter $C_{(k)}$ is a notation of the sequence $0\leq c_{(k)}(1)\leq\cdots\leq c_{(k)}(m_k)\leq n$,
according to our conventions on simplicial sets.

Explicitly, the interval $\{0,\ldots,n\}$ is cut at positions specified by indices
$$0 = n_0\leq n_1\leq\cdots\leq n_{r+d-1}\leq n_{r+d} = n$$
as represented by the following diagram
$$\xymatrix@M=0pt@C=12mm{ \ar@{|-|}[r]^<{0} & \ar@{|-|}[r]^<{n_1} &
\ar@{.}[r]^<{n_2} & \ar@{.}[r] & \ar@{|-|}[r]^>{n} & \\ }$$
The sequence $C_{(k)}$ is the concatenation of the intervals $[n_{i-1},n_i]$
such that $u(i) = k$.
Hence, if $u(i_1),\ldots,u(i_e)$ denote the occurences of $k$
in the sequence $(u(1),\ldots,u(r+d))$,
then the element $\Delta(C_{(k)})$ is represented by the simplex:
$$\Delta(n_{i_1-1}\,\seg\,n_{i_1},n_{i_2-1}\,\seg\,n_{i_2},\ldots,n_{i_e-1}\,\seg\,n_{i_e}).$$
This simplex is non-degenerate if $n_{i_{*-1}}<n_{i_*-1}$,
for all $* = 1,\ldots,e-1$.
In paragraph \ref{IntervalLength},
we identify the dimension of this simplex to the sum of the lengths
of the intervals $[n_{i_*-1},n_{i_*}]$:
$$\eqalign{ m_k &
= \dim\Delta(n_{i_1-1}\,\seg\,n_{i_1},n_{i_2-1}\,\seg\,n_{i_2},\ldots,n_{i_e-1}\,\seg\,n_{i_e}) \cr
& = \length_u(n_{i_1-1}\,\seg\,n_{i_1}) + \length_u(n_{i_2-1}\,\seg\,n_{i_2})
+ \cdots + \length_u(n_{i_e-1}\,\seg\,n_{i_e}). \cr }$$

\paragraph{\it Example}\label{ExIntervalCut}
Here is a convenient graphical representation which we adopt in the article.
The intervals $[n_{i-1},n_i]$ are arranged on $r$ lines
in a diagram which has $r+d$ columns delimited by the variables $n_i$.
The interval $[n_{i-1},n_i]$ is put on the line $u(i)$
in the $i$th column of the diagram.
As an example, the surjection, which is specified by the sequence $(3,1,2,3,2)$,
gives the interval cut diagram
$$
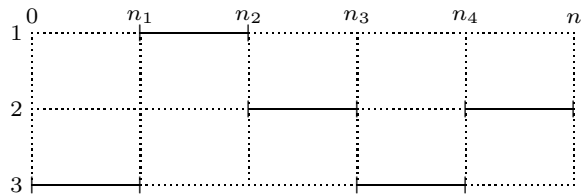
\xymatrix@M=0pt@C=12mm{\ar@{.}[d]_<{1}\ar@{.}[r]^<{0} &
\ar@{.}[d]\ar@{|-|}[r]^<{n_1} &
\ar@{.}[d]\ar@{.}[r]^<{n_2} &
\ar@{.}[d]\ar@{.}[r]^<{n_3} &
\ar@{.}[d]\ar@{.}[r]^<{n_4}^>{n} &
\ar@{.}[d] \\
\ar@{.}[d]_<{2}_>{3}\ar@{.}[r] &
\ar@{.}[d]\ar@{.}[r] &
\ar@{.}[d]\ar@{|-|}[r] &
\ar@{.}[d]\ar@{.}[r] &
\ar@{.}[d]\ar@{|-|}[r] &
\ar@{.}[d] \\ 
\ar@{|-|}[r] & \ar@{.}[r] & \ar@{.}[r] & \ar@{|-|}[r] & \ar@{.}[r] & \\ }$$
and is associated to the tensor
$$\Delta(n_1\,\seg\,n_2)
\otimes\Delta(n_2\,\seg\,n_3,n_4\,\seg\,n)
\otimes\Delta(0\,\seg\,n_1,n_3\,\seg\,n_4).$$
The dimensions of these simplices are respectively $m_1 = (n_2-n_1)$,
$m_2 = (n_3-n_2+1)+(n-n_4)$ and $m_3 = (n_1+1)+(n_4-n_3)$.
Equivalently, the diagram is a representation of the morphisms
$$\matrix{ & \{0,\ldots,m_1\} &
\rto{\simeq}{} & \{n_1,\ldots,n_2\}\hfill &
\,\hookrightarrow\, & \{0,\ldots,n\}, \cr
& \{0,\ldots,m_2\} &
\rto{\simeq}{} & \{n_2,\ldots,n_3,n_4,\ldots,n\}\hfill &
\,\hookrightarrow\, & \{0,\ldots,n\}, \cr
\hbox{and}\qquad & \{0,\ldots,m_3\} &
\rto{\simeq}{} & \{0,\ldots,n_1,n_3,\ldots,n_4\}\hfill &
\,\hookrightarrow\, & \{0,\ldots,n\}, \cr }$$
which are denoted by $c_{(k)}: [m_k]\,\rightarrow\,[n]$ in the definition.

\paragraph{\it On the interval length determined by an interval cut}\label{IntervalLength}
In the interval cut associated to a surjection,
there are {\it inner} and {\it final} intervals.
The interval $[n_{i-1},n_i]$ is final,
if $u(i)$ is the last occurrence of a value in the surjection.
The interval $[n_{i-1},n_i]$ is inner,
if $u(i)$ is a caesura of the surjection.
Consider again the simplex:
$$\Delta(n_{i_1-1}\,\seg\,n_{i_1},n_{i_2-1}\,\seg\,n_{i_2},\ldots,n_{i_e-1}\,\seg\,n_{i_e})$$
which is associated to the occurences of $k$ in the surjection $u$.
The interval $[n_{i_x-1},n_{i_x}]$ is inner for $x = 1,\ldots,e-1$
and final for $x = e$.

The length of an interval $[n_{i-1},n_i]$
is defined by $\length_u [n_{i-1},n_i] = n_i-n_{i-1}$ for a final interval
and by $\length_u [n_{i-1},n_i] = n_i-n_{i-1}+1$ for an inner interval.
Hence, in the context above, we have
$$\eqalign{ \length_u(n_{i_x-1}\,\seg\,n_{i_x}) &
= n_{i_x}-n_{i_x-1}+1,\qquad\hbox{for}\ x = 1,\ldots,e-1, \cr
\hbox{and}\qquad\length_u(n_{i_x-1}\,\seg\,n_{i_x}) &
= n_{i_x}-n_{i_x-1},\qquad\hbox{for}\ x= e. \cr }$$
Therefore, we have the formula:
$$\eqalign{ & \dim\Delta(n_{i_1-1}\,\seg\,n_{i_1},n_{i_2-1}\,\seg\,n_{i_2},
\ldots,n_{i_e-1}\,\seg\,n_{i_e}) \cr
= & \length_u(n_{i_1-1}\,\seg\,n_{i_1}) + \length_u(n_{i_2-1}\,\seg\,n_{i_2})
+ \cdots + \length_u(n_{i_e-1}\,\seg\,n_{i_e}). \cr }$$

\paragraph{\it On the signs associated to an interval cut}\label{SignIntervalCut}
The interval cut defined in paragraph \ref{IntervalCut} has an associated
permutation sign and an associated position sign.
The permutation sign depends on the relative position
of the intervals $[n_{i-1},n_i]$.
The position sign depends explicitly on the absolute position
of the inner intervals  $[n_{i-1},n_i]$.
The total sign associated to an interval cut is the product
of the permutation sign and the position sign.

The permutation sign is determined as follows.
Explicitly, we consider the shuffle which takes $(u(1),\ldots,u(r+d))$
to the ordered sequence $(1,\ldots,1,2,\ldots,2,\ldots,r,\ldots,r)$.
The permutation of the associated intervals produces a sign
according to the permutation rule in dg-calculus
with the lengths of the intervals as degrees.
This sign is the permutation sign associated to the interval cut.

The position signs are provided by the inner intervals.
Explicitly, if the interval $[n_{i-1},n_i]$ is inner,
then this interval has $n_i$ as an associated position sign-exponent.

\paragraph{\it Example}
As an example, take the interval cut introduced in paragraph \ref{ExIntervalCut},
which is associated to the sequence $u = (3,1,2,3,2)$:
$$
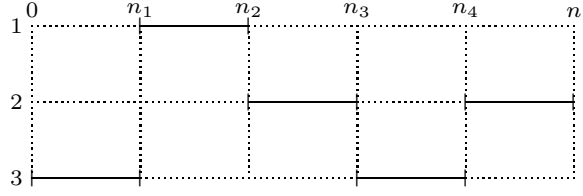
\xymatrix@M=0pt@C=12mm{\ar@{.}[d]_<{1}\ar@{.}[r]^<{0} &
\ar@{.}[d]\ar@{|-|}[r]^<{n_1} &
\ar@{.}[d]\ar@{.}[r]^<{n_2} &
\ar@{.}[d]\ar@{.}[r]^<{n_3} &
\ar@{.}[d]\ar@{.}[r]^<{n_4}^>{n} &
\ar@{.}[d] \\
\ar@{.}[d]_<{2}_>{3}\ar@{.}[r] &
\ar@{.}[d]\ar@{.}[r] &
\ar@{.}[d]\ar@{|-|}[r] &
\ar@{.}[d]\ar@{.}[r] &
\ar@{.}[d]\ar@{|-|}[r] &
\ar@{.}[d] \\ 
\ar@{|-|}[r] & \ar@{.}[r] & \ar@{.}[r] & \ar@{|-|}[r] & \ar@{.}[r] & \\ }$$
In this case, the final intervals are $[n_1,n_2]$, $[n_3,n_4]$ and $[n_4,n]$.
These interval have length:
$$\eqalign{ & \length_u(n_1\,\seg\,n_2) = n_2-n_1,\qquad\length_u(n_3\,\seg\,n_4) = n_4-n_3 \cr
& \hbox{and}\qquad\length_u(n_4\,\seg\,n) = n-n_4. \cr }$$
The inner intervals are $[0,n_1]$ and $[n_2,n_3]$ and have length:
$$\length_u(0\,\seg\,n_1) = n_1+1\qquad\hbox{and}\qquad\length_u(n_2\,\seg\,n_3) = n_3-n_2+1.$$
The associated position sign is given by the exponent $\tau = n_1+n_3$.
The permutation sign is determined by the shuffle:
$$[n_1,n_2],[n_2,n_3],[n_4,n],[0,n_1],[n_3,n_4]
\mapsto [0,n_1],[n_1,n_2],[n_2,n_3],[n_3,n_4],[n_4,n]$$
In diagrams:
$$\vcenter
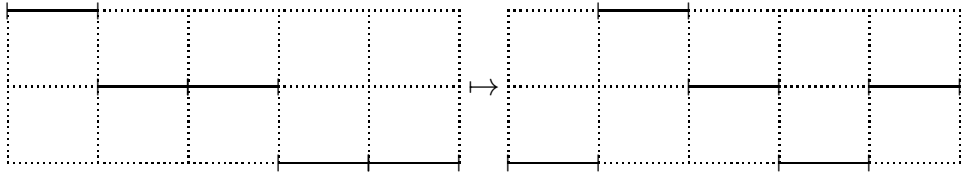
{\xymatrix@M=0pt@C=10mm{\ar@{.}[d]\ar@{|-|}[r] &
\ar@{.}[d]\ar@{.}[r] &
\ar@{.}[d]\ar@{.}[r] &
\ar@{.}[d]\ar@{.}[r] &
\ar@{.}[d]\ar@{.}[r] &
\ar@{.}[d] \\
\ar@{.}[d]\ar@{.}[r] &
\ar@{.}[d]\ar@{|-|}[r] &
\ar@{.}[d]\ar@{|-|}[r] &
\ar@{.}[d]\ar@{.}[r] &
\ar@{.}[d]\ar@{.}[r] &
\ar@{.}[d] \\ 
\ar@{.}[r] & \ar@{.}[r] & \ar@{.}[r] & \ar@{|-|}[r] & \ar@{|-|}[r] & \\ }}
\mapsto\vcenter{\xymatrix@M=0pt@C=10mm{\ar@{.}[d]\ar@{.}[r] &
\ar@{.}[d]\ar@{|-|}[r] &
\ar@{.}[d]\ar@{.}[r] &
\ar@{.}[d]\ar@{.}[r] &
\ar@{.}[d]\ar@{.}[r] &
\ar@{.}[d] \\
\ar@{.}[d]\ar@{.}[r] &
\ar@{.}[d]\ar@{.}[r] &
\ar@{.}[d]\ar@{|-|}[r] &
\ar@{.}[d]\ar@{.}[r] &
\ar@{.}[d]\ar@{|-|}[r] &
\ar@{.}[d] \\ 
\ar@{|-|}[r] & \ar@{.}[r] & \ar@{.}[r] & \ar@{|-|}[r] & \ar@{.}[r] & \\ }}$$
Hence, the permutation sign-exponent associated to the interval cut is the sum
$\sigma = (n_1+1)(n_2-n_1) + (n_1+1)(n_3-n_2+1) + (n_1+1)(n-n_4) + (n_4-n_3)(n-n_4)$.

\paragraph{\it The interval cut operation associated to a surjection}
The interval cut operation $\AW(u)\in\Z(r)$ associated to $u\in\X(r)_d$
maps the simplex $\Delta(0,\ldots,n)\in\Delta^n$
to the sum
$$\AW(u)(\Delta(0,\ldots,n)) = \sum\pm\Delta(C_{(1)})\otimes\cdots\otimes\Delta(C_{(r)})$$
over the indices $0 = n_0\leq n_1\leq\cdots\leq n_{r+d-1}\leq n_{r+d} = n$,
where $\Delta(C_{(1)}),\ldots,\Delta(C_{(r)})$ are the interval cut simplices
defined in paragraph \ref{IntervalCut}.
The sign associated to a term is the product of the permutation and position signs
defined in paragraph \ref{SignIntervalCut}.

In general, the image of a simplex $x\in X_n$ is determined by the formula:
$$\AW(u)(x) = \sum\pm x(C_{(1)})\otimes\cdots\otimes x(C_{(r)}).$$
As an example, if $u = (3,1,2,3,2)$, then we have
$$\AW(u)(x) = \sum\pm x(n_1\,\seg\,n_2)\otimes x(n_2\,\seg\,n_3,n_4\,\seg\,n)
\otimes x(0\,\seg\,n_1,n_3\,\seg\,n_4)$$
for any $n$-dimensional simplex $x\in X_n$.
Equivalently, we have:
$$\AW(u)(x) = \sum\pm c_{(1)}^*(x)\otimes\cdots\otimes c_{(r)}^*(x),$$
where $c_{(1)},\ldots,c_{(r)}$ are the morphisms in the simplicial category
represented by the sequences $C_{(1)},\ldots,C_{(r)}$.

\paragraph{\sc Lemma}
{\it The map $u\mapsto\AW(u)$ defines an operad morphism $\AW: \X\,\rightarrow\,\Z$.}

\smallskip
We omit the proof of this lemma which reduces to tedious verifications.
We refer to J. McClure and J. Smith ({\it cf}. [\cite{MCSSeq}]) for details.
In our conventions,
it is straighforward to verify that our signs are coherent
using the invariance properties of the length
in an interval cut diagram.

\paragraph{\it On the Alexander-Whitney diagonal}
The operation $\AW(1,2): N_*(X)\,\rightarrow\,N_*(X)^{\otimes 2}$
is identified with the Alexander-Whitney diagonal.
We have clearly:
$$\AW(1,2)(x) = \sum_{0\leq i\leq n} x(0\,\seg\,i)\otimes x(i\,\seg\,n).$$
Similarly, the operations
$$\AW(\theta_d): N_*(X)\,\rightarrow\,N_*(X)^{\otimes 2},$$
where $\theta_d = (1,2,1,2,\ldots)\in\X(2)_d$,
are identified with the higher Alexander-Whitney diagonals,
which induce the cup-$d$-products on cochains
({\it cf}. [\cite{S}]).
We have explicitly:
$$\matrix{ \AW(1,2,1)(x)\hfill &
\displaystyle = \sum\nolimits_{i<j} \pm x(0\,\seg\,i,j\,\seg\,n)\otimes x(i\,\seg\,j),\hfill\cr
\AW(1,2,1,2)(x)\hfill &
\displaystyle = \sum\nolimits_{i<j<k} \pm x(0\,\seg\,i,j\,\seg\,k)
\otimes x(i\,\seg\,j,k\,\seg\,n),\hfill\cr }$$
and so on.
The sign-exponent is $(n-j)(j-i)+i$ for $\AW(1,2,1)$ and $(k-j)(j-i+1)+i+j$ for $\AW(1,2,1,2)$.

It should be clear that the element $\mu_d = (\id,\tau,\id,\tau,\ldots)\in\E(r)_d$
introduced in theorem \ref{BarrattEcclesChainOperations} verifies $\TR(\mu_d) = \theta_d$.
We conclude that this element gives a representative of the cup-$d$-product
as claimed in theorem \ref{BarrattEcclesChainOperations}.

\section{On closed model structures}

\subsection{The closed model structure}

\smallskip
The purpose of this section is to prove the following theorem:

\paragraph{\sc Theorem}\label{Models}
{\it The $\E$-algebras form a closed model category.
A weak equivalence (respectively, a fibration)
is an algebra morphism
which is a weak equivalence (respectively, a fibration)
in the category of dg-modules.}

\smallskip
In fact, our arguments extend to a wide class of operads.
In general, an operad $\P$ has a canonical $\Sigma_*$-projective resolution
which is provided by the tensor product with the Barratt-Eccles operad.
More explicitly, we consider the operad such that $(\E\otimes\P)(r) = \E(r)\otimes\P(r)$.
We have a canonical augmentation $\E\otimes\P\,\rightarrow\,\P$
which is induced by the augmentation $\E\,\rightarrow\,\C$
of the Barratt-Eccles operad. Our arguments extend to the following situation:

\paragraph{\sc Assumption:}\label{Assumptn}
{\it The augmentation morphism $\E\otimes\P\,\rightarrow\,\P$ has a section $\P\,\rightarrow\,\E\otimes\P$
in the category of dg-operads.}

\smallskip
To be explicit:

\paragraph{\sc Theorem}\label{CMOperads}
{\it If the operad $\P$ satisfies assumption \ref{Assumptn},
then the $\P$-algebras form a closed model category.
A weak equivalence (respectively, a fibration)
is an algebra morphism
which is a weak equivalence (respectively, a fibration)
in the category of dg-modules.}

\smallskip
The assumption is satisfied by the Barratt-Eccles operad $\E$.
In this case, a section $\E\,\rightarrow\,\E\otimes\E$ is given by the Alexander-Whitney diagonal
(see paragraph \ref{AWDiag}).
More generally, we have:

\paragraph{\sc Fact:}
{\it Assumption \ref{Assumptn} is satisfied by the classical $\Sigma_*$-projective resolution
$\P = \E\otimes\Q$ of an operad $\Q$.
The section $\P\,\rightarrow\,\E\otimes\P$
is provided by the diagonal of the Barratt-Eccles operad
$\E\otimes\Q\,\rightarrow\,\E\otimes\E\otimes\Q$.}

\smallskip
The assumption is also satisfied in the following situation:

\paragraph{\sc Fact:} {\it Assumption \ref{Assumptn} is satisfied by cofibrant operads $\P$.
In this case, the existence of the section $\P\,\rightarrow\,\E\otimes\P$
follows from the left lifting property of operad cofibrations
({\it cf}. V. Hinich, [\cite{H}]).}

\smallskip
By definition, a morphism of dg-modules is a fibration if it is surjective (in all degrees)
and a weak equivalence if it induces an isomorphism in homology.
The cofibrations are the morphisms which have the left lifting property with respect to acyclic fibrations.

We deduce theorem \ref{CMOperads} from standard arguments,
which go back to Quillen (cf. [\cite{Q}]),
and we omit most of the demonstration
(we refer also to [V. Hinich, \cite{H}] and [M. Livernet, \cite{L}]).
We prove only the following lemma which is the crucial point
of the proof of theorem \ref{CMOperads}:

\paragraph{\sc Lemma}\label{LLP}
{\it If assumption \ref{Assumptn} is satisfied, then a $\P$-algebra morphisms
which has the left lifting property with respect to fibrations is a weak equivalence.}

\smallskip
We would like to mention that lemma \ref{LLP} (and, accordingly, theorem \ref{CMOperads})
does not hold for any operad.
In general,
we have just the following weaker property
({\it cf}. [V. Hinich, \cite{H}], [M. Mandell, \cite{M}]):

\paragraph{\sc Fact:} {\it We assume that $\P$ is a $\Sigma_*$-projective operad.
As long as $A$ is a cofibrant algebra,
a $\P$-algebra morphism $f: A\,\rightarrow\,B$
which has the left lifting property with respect to fibrations
is a weak equivalence.}

\smallskip
First, we observe that, in the situation of assumption \ref{Assumptn},
any $\P$-algebra has a canonical path object.
Then, we are able to deduce lemma \ref{LLP} from an argument of D. Quillen
(as in characteristic $0$, {\it cf}. [\cite{L}]).

\paragraph{\it On path objects}
Let $A$ be a $\P$-algebra.
Recall that a path object of $A$ is a $\P$-algebra $\tilde{A}$ together with morphisms
$$\xymatrix{ A\ar[r]^{\sim}_{s_0} & \tilde{A}\ar@<0.5ex>[r]^{d_0}\ar@<-0.5ex>[r]_{d_1} & A }$$
such that $s_0$ is a weak equivalence, $(d_0,d_1): \tilde{A}\,\rightarrow\,A\times A$ is a fibration
and $d_0 s_0 = d_1 s_0 = 1_A$.
If the operad $\P$ satisfies assumption \ref{Assumptn},
then, for any $\P$-algebra $A$, the tensor products
$$\xymatrix{ N^*(\Delta^0)\otimes A\ar[r]^{\sim}_{s_0\otimes A} &
N^*(\Delta^1)\otimes A\ar@<0.5ex>[r]^{d_0\otimes A}\ar@<-0.5ex>[r]_{d_1\otimes A} &
N^*(\Delta^0)\otimes A }$$
give a path-object $\tilde{A} = N^*(\Delta^1)\otimes A$.

In fact, the tensor product $K\otimes A$ of a $\P$-algebra $A$
with an $\E$-algebra $K$ forms an $\E\otimes\P$-algebra.
If there is an operad morphism $\P\,\rightarrow\,\E\otimes\P$,
then the tensor product $K\otimes A$ is a $\P$-algebra
by restriction of structure.
Furthermore, the $\E$-algebra $N^*(\Delta^0)$ is identified with the ground ring $\F$,
which is equipped with the $\E$-algebra structure obtained by restriction
through the operad augmentation $\E\,\rightarrow\,\C$.
Consequently, the tensor product $N^*(\Delta^0)\otimes A$ is identified with the dg-module $A$
together with the $\E\otimes\P$-algebra structure obtained by restriction
through the augmentation $\E\otimes\P\,\rightarrow\,\P$.
Therefore, if the morphism $\P\,\rightarrow\,\E\otimes\P$ is a section
of the augmentation $\E\otimes\P\,\rightarrow\,\P$,
then the dg-module $N^*(\Delta^0)\otimes A$ is identified with $A$
as a $\P$-algebra.
We conclude that the tensor product $N^*(\Delta^1)\otimes A$
is a path object for $A$.
Let us mention that the tensor product $N^*(\Delta^1)\otimes A$
is identified with the classical cylinder construction
in the category of dg-modules.

\smallskip
To recapitulate:

\paragraph{\sc Proposition}
{\it If assumption \ref{Assumptn} is satisfied,
then a $\P$-algebra $A$ has a canonical path object
which is represented by the tensor product $N^*(\Delta^1)\otimes A$.}

\smallskip
Lemma \ref{LLP} is a consequence of the following classical observation
({\it cf}. [\cite{L}, Chapter II, Lemma 2.5.4], [\cite{Q}, Section I.5, Lemma 1]):

\paragraph{\sc Observation:}
{\it The morphisms which have the left lifting property with respect to fibrations
are strong deformations retracts.}

\smallskip
Just recall that a morphism $f: A\,\rightarrow\,F$ is a strong deformation retract
if there is a morphism $r: F\,\rightarrow\,A$ such that $r f = 1_F$
and a morphism $h: F\,\rightarrow\,N^*(\Delta^1)\otimes F$
such that $d_0 h = f r$, $d_1 h = 1_F$ and $h f = s_0 f$.
Clearly, these identities imply that $f: A\,\rightarrow\,F$ is a weak equivalence.

\subsection{On spheres, cones and suspensions}

\paragraph{\it On spheres}
We let $S^n$ denote the standard simplicial model of the $n$-dimensional sphere:
$$S^n = \Delta^n/\cup_{i=0}^n\Delta^n(0,\ldots,i-1,i+1,\ldots,n).$$
This simplicial set has just a non-degenerate simplex $\Delta(0,\ldots,n)\in S^n$
in dimension $n$
and a base point $\pt\in S^n$ in dimension $0$.
We let $S(n)$ denote the reduced normalized cochain complex
$S(n) = \tilde{N}^*(S^n)$ associated to $S^n$.
In fact, the reduced chain complex $\tilde{N}_*(S^n)$ is concentrated in degree $n$
and is generated by $\Delta(0,\ldots,n)$.
We let ${\bf e}^n\in\bar{N}^*(S^n)$ denote the dual element.
Hence, we have:
$$S(n) = \F\cdot{\bf e}^n.$$
If $n = 1$, then, for simplicity, we write ${\bf e} = {\bf e}^1$.

The $\E$-algebra product
$$\E(r)_d\otimes(S(n)^{\otimes r})^*\,\rightarrow\,S(n)^{*-d}$$
has a non-trivial component in degree $*=n r$ and $d=n(r-1)$.
One purpose of this section is to make this component explicit.

\paragraph{\it On the pointed interval}
Similarly, we let $C(1)$ denote the reduced normalized cochain complex
of the standard simplicial interval $\Delta^1$
equipped with the base point $* = \Delta(0)$.
In fact, the reduced chain complex $\tilde{N}_*(\Delta^1)$
is generated by $\Delta(0,1)\in\tilde{N}_1(C^n)$ and $\Delta(1)\in\tilde{N}_0(C^n)$.
Therefore, if ${\bf e}\in\bar{N}^1(\Delta^1)$
and ${\bf c}\in\bar{N}^0(\Delta^1)$ denote the dual elements,
then we have:
$$C(1) = \F\cdot{\bf e}\oplus\F\cdot{\bf c},$$
together with the differential $\delta({\bf c}) = {\bf e}$.
We make also the $\E$-algebra products in $C(1)$ explicit.

To conclude this paragraph, we have the canonical exact sequence:
$$0\,\rightarrow\,S(1)\,\rightarrow\,C(1)\,\rightarrow\,S(0)\,\rightarrow\,0.$$
More explicitly, the circle $S(1)$ is isomorphic to the sub-dg-module
$\F\cdot{\bf e}\,\hookrightarrow\,C(1)$
and the zero-sphere $S(0)$ is isomorphic to the quotient dg-module
$C(1)\,\rightarrow\,\F\cdot{\bf c}$.

\paragraph{\it Certain fundamental cochains}
We introduce a cochain $\epsilon_s: \E(r)_*\,\rightarrow\,\F$
in order to make the algebra structure of $S(1)$ and $C(1)$ explicit.
Consider a $d$-dimensional simplex $(w_0,\ldots,w_d)\in\E(r)_d$.
If the sequence $(w_0(1),\ldots,w_d(1))$ does not form a permutation of $(1,\ldots,s)$,
then we set $\epsilon_s(w_0,\ldots,w_d) = 0$.
(In particular,
we assume that the cochain $\epsilon_s: \E(r)_*\,\rightarrow\,\F$
vanishes in degree $d\not=s-1$.)
Otherwise,
we set $\epsilon_s(w_0,\ldots,w_{s-1}) = \pm 1$
according to the signature of $(w_0(1),\ldots,w_{s-1}(1))$.
Explicitly:
$$\epsilon_s(w_0,\ldots,w_{s-1}) = \sgn(w_0(1),\ldots,w_{s-1}(1)).$$
By convention,
the $0$-cochain $\epsilon_0: \E(r)_*\,\rightarrow\,\F$ is the classical augmentation
of the bar complex $\E(r)_*$.

\paragraph{\sc Theorem}\label{SphereModels}
{\it 1) The evaluation product of the circle algebra $S(1)$
is given by the following formula:
$$w({\bf e}_{(1)},\ldots,{\bf e}_{(r)}) = (-1)^\sigma\cdot\epsilon_r(w)\cdot{\bf e},\leqno{i)}$$
where $\sigma = r(r-1)/2$.

2) The evaluation product of the cone algebra $C(1)$ is characterized
by the following equations:
$$\leqalignno{ & w({\bf e}_{(1)},\ldots,{\bf e}_{(r)})
= (-1)^\sigma\cdot\epsilon_r(w)\cdot{\bf e}, & i) \cr
& w({\bf e}_{(1)},\ldots,{\bf e}_{(s)},{\bf c}_{(s+1)},\ldots,{\bf c}_{(r)})
= (-1)^\sigma\cdot\epsilon_s(w)\cdot{\bf e} & ii) \cr
& \hbox{and}\qquad w({\bf c}_{(1)},\ldots,{\bf c}_{(r)})
= (-1)^\sigma\cdot\epsilon_0(w)\cdot{\bf c} & iii) \cr }$$
(where, whenever it makes sense, $0<s<r$).
The subscripts in these formulas specify the places of the copies of ${\bf e}$ and ${\bf c}$.
In fact, equation i) is the instance $s = r$ of equation ii).
In all cases, the sign-exponent is given by $\sigma = s(s-1)/2$
where $s = r$ for equation i) and $s = 0$ for equation iii).}

\smallskip
We have morphisms of dg-modules
$$\matrix{ & \tilde{N}^*(S^1\wedge X) &
\rightarrow & \tilde{N}^*(S^1)\otimes\tilde{N}^*(X) &
= & \Sigma^*\tilde{N}^*(X) \cr
\hbox{and}\qquad & \tilde{N}^*(\Delta^1\wedge X) &
\rightarrow & \tilde{N}^*(\Delta^1)\otimes\tilde{N}^*(X) &
= & C^*\tilde{N}^*(X) \cr }$$
given by the classical Eilenberg-Zilber equivalence.
But, these morphisms are not compatible with $\E$-algebra structures in general.

Nevertheless,
it is possible to generalize the formulas of theorem \ref{SphereModels}
to the $n$-dimensional spheres $S(n)$.
In fact, we obtain the following result, which we mention as a remark:

\paragraph{\sc Proposition}\label{HigherDimSpheres}
{\it The canonical isomorphim $S(n)\,\rightarrow\,S(1)^{\otimes n}$
which identifies the element ${\bf e}^n\in S(n)$
with the tensor ${\bf e}^{\otimes n}\in S(1)^{\otimes n}$
is an isomorphism of $\E$-algebras.}

\paragraph{\it On cones and suspensions}\label{Cone}
If $V\in\dg\Mod_\F$ is an (upper-graded) dg-module,
then $C^* V\in\dg\Mod_\F$ denotes the classical cone associated to $V$
in the category of dg-modules.
We have in fact $C^* V = C(1)\otimes V$.
Accordingly:
$$(C^* V)^d = {\bf c}\otimes V^d\oplus{\bf e}\otimes V^{d-1}$$
and the differential of $C^* V$ is given by
$$\delta({\bf c}\otimes x+{\bf e}\otimes y)
= {\bf c}\otimes\delta(x) + {\bf e}\otimes(x-\delta(y)).$$
Similarly, the suspension of $V$ in the category of dg-modules,
denoted by $\Sigma^* V\in\dg\Mod_\F$,
is identified with the tensor product $S(1)\otimes V$.
Equivalently:
$$(\Sigma^* V)^d = {\bf e}\otimes V^{d-1}.$$
We have a canonical short exact sequence:
$$0\,\rightarrow\,\Sigma^* V\,\rightarrow\,C^* V\,\rightarrow\,V\,\rightarrow\,0.$$
If $A$ is an $\E$-algebra,
then $C^* A = C(1)\otimes A$ is a tensor product of $\E$-algebras.
Therefore,
the cone $C^* A$ has the structure of an $\E$-algebra.
Similarly,
the suspension $\Sigma^* A = S(1)\otimes A$ has a natural $\E$-algebra structure.
We obtain readily:

\paragraph{\sc Proposition}
{\it Let $A$ be an $\E$-algebra.
The cone $C^* A = {\bf e}\otimes A\oplus{\bf c}\otimes A$
is equipped with a natural $\E$-algebra structure.
Furthermore, the suspension $\Sigma^* A = {\bf e}\otimes A\hookrightarrow C^* A$
is a subalgebra of $C^* A$
and the canonical surjection $C^* A\rightarrow {\bf c}\otimes A = A$
is a morphism of $\E$-algebras.
In $C^* A$, the products are given by the following formulas:
$$\eqalign{ & w({\bf e}\otimes a_{(1)},\ldots,{\bf e}\otimes a_{(r)})
= \pm(-1)^\sigma\cdot{\bf e}\otimes(\epsilon_r\cap w)(a_{(1)},\ldots,a_{(r)}), \cr
& w({\bf e}\otimes a_{(1)},\ldots,{\bf e}\otimes a_{(s)},
{\bf c}\otimes a_{(s+1)},\ldots,{\bf c}\otimes a_{(r)})
= \pm(-1)^\sigma\cdot{\bf e}\otimes(\epsilon_s\cap w)(a_{(1)},\ldots,a_{(r)}), \cr
& \hbox{and}\qquad w({\bf c}\otimes a_{(1)},\ldots,{\bf c}\otimes a_{(r)})
= \pm(-1)^\sigma\cdot{\bf c}\otimes w(a_{(1)},\ldots,a_{(r)}), \cr }$$
where $0<s<r$.
In fact, the first formula is the instance $s = r$ of the second formula.
But, the formula for $s = 0$ differs from the general case.

In all cases, we have $\sigma = s(s-1)/2$ (as in theorem \ref{SphereModels}).
The unspecified signs are determined by the commutation of the elements ${\bf e}$
with the factors $w,a_{(1)},\ldots,a_{(r)}$.
Accordingly, these signs are determined by the sign exponent
$\deg(w)\cdot s+\deg(a_1)\cdot(s-1)+\cdots+\deg(a_{s-1})\cdot 1$.}

\smallskip
In our context, the cap product of a cochain $\phi: \E(r)_d\,\rightarrow\,\F$
with a chain $w = (w_0,\ldots,w_n)\in\E(r)_n$
is represented by the element
$$\phi\cap w = \phi(w_0,\ldots,w_d)\cdot(w_d,\ldots,w_n)\in\E(r)_{n-d}.$$

\paragraph{\it On the operad suspension}\label{OperadSuspension}
There is an operad $\Lambda^*\E$, associated to $\E$,
whose algebras are suspensions $\Sigma^* A$ of $\E$-algebras $A$.
We have in fact:
$$\Lambda^*\E(r)_d = \sgn(r)\otimes\E(r)_{d-r+1}$$
where $\sgn(r)$ denotes the signature representation ({\it cf}. [\cite{GeJ}]).
If the letter $w$ denotes an element of $\E(r)_d$,
then $\Lambda^* w$ is the associated operation in $\Lambda^*\E(r)_{d+r-1}$.
The suspension $\Sigma^* A$ of an $\E$-algebra $A$ is equipped with the evaluation product
$$\Lambda^*\E(r)\otimes(\Sigma^* A)^{\otimes r}\,\rightarrow\,(\Sigma^* A)$$
such that
$$(\Lambda^* w)({\bf e}\otimes a_1,\ldots,{\bf e}\otimes a_r)
= \pm {\bf e}\otimes w(a_1,\ldots,a_r).$$
The sign is produced by the commutation of the elements ${\bf e}$
with the factors $w,a_1,\ldots,a_r$ and is determined
by the parity of the sign-exponent
$\deg(w)\cdot r+\deg(a_1)\cdot(r-1)+\cdots+\deg(a_{r-1})\cdot 1+\deg(a_r)\cdot 0$.

The previous proposition has the following consequence:

\paragraph{\sc Proposition}\label{SuspensionMorphism}
{\it The cap products $\epsilon_r\cap-: \E(r)_*\,\rightarrow\,\E(r)_{*-r+1}$, $r>0$,
define an operad morphism $\epsilon_*\cap-: \E\,\rightarrow\,\Lambda^*\E$.}

\paragraph{\it Remark}
The operads $\X$ and $\Z$ are also equipped with such morphisms
$\X\,\rightarrow\,\Lambda^*\X$ and $\Z\,\rightarrow\,\Lambda^*\Z$.
Furthermore, we have a commutative diagram:
$$\xymatrix{ \E\ar[r]^{\TR}\ar[d] & \X\ar[r]^{\AW}\ar@{-->}[d] & \Z\ar@{-->}[d] \\
\Lambda^*\E\ar@{->}[r]_{\Lambda^*\TR} & \Lambda^*\X\ar@{->}[r]_{\Lambda^*\AW} & \Lambda^*\Z \\ }$$

The morphism $\X(r)_*\,\rightarrow\,\X(r)_{*-r+1}$
is very similar to $\epsilon_r\cap-: \E(r)_*\,\rightarrow\,\E(r)_{*-r+1}$.
Precisely, given $u\in\X(r)_d$, we let $\epsilon_r\cap u\in\X(r)_{d-r+1}$
be the sequence such that
$$\epsilon_r\cap u = \sgn(u(1),\ldots,u(r))\cdot(u(r),u(r+1),\ldots,u(r+d)),$$
if $(u(1),\ldots,u(r))$ is a permutation of $(1,\ldots,r)$,
and $\epsilon_r\cap u = 0$, otherwise.
This cap product is an operad morphism
$\epsilon_*\cap-: \X\,\rightarrow\,\Lambda^*\X$
and makes commute the left-hand square of the diagram above.

The morphism $\Z(r)_*\,\rightarrow\,\Z(r)_{*-r+1}$ is due to V. Smirnov
({\it cf}. [\cite{SmCLoop}], [\cite{SmHLoop}]).
To be explicit,
we let $\sigma: N_*(\Delta^n)\,\rightarrow\,N_{*-1}(\Delta^{n-1})$
be the dg-morphism such that:
$$\sigma(\Delta(i_0,i_1,\ldots,i_n)) = \left\{\matrix{ \Delta(i_1-1,\ldots,i_n-1),\hfill & \ \hbox{if}\ i_0 = 0,\hfill \cr
0,\hfill & \ \hbox{otherwise}.\hfill \cr }\right.$$
The tensor powers
$$\sigma^{\otimes r}: (N_*(\Delta^n)^{\otimes r})_d\,\rightarrow\,(N_*(\Delta^{n-1})^{\otimes r})_{d-r}$$
are equivalent to a morphism
$$\Z(r)_{d-n}\,\rightarrow\,\Z(r)_{d-r+1-n}.$$
In fact, this morphism of dg-modules is a morphism of operads $\Z\,\rightarrow\,\Lambda^*\Z$.
It is also straightforward to prove that it makes commute
the right-hand square of the diagram.

\subsection{Some proofs}

\smallskip
In this section, we determine the structure of the chain complexes
$$\tilde{N}_*(X)\ =\ \tilde{N}_*(S^0),\ \tilde{N}_*(\Delta^1)\ \hbox{and}\ \tilde{N}_*(S^1)$$
as coalgebras over the surjection operad $\X$.
The algebra structure of $\tilde{N}^*(X)$, as specified in theorem \ref{SphereModels},
follows from these calculations and from our definitions.

\smallskip
To begin with, we record the following fact:

\paragraph{\sc Fact:}
{\it There are morphisms of pointed simplicial sets
$S^0\,\rightarrow\,\Delta^1\,\rightarrow\,S^1$
which induce morphisms of coalgebras:
$$\tilde{N}_*(S^0)\,\rightarrow\,\tilde{N}_*(\Delta^1)\,\rightarrow\,\tilde{N}_*(S^1).$$
We have explicitly:
$$\tilde{N}_*(S^0)
= \F\cdot\Delta(0),\qquad\tilde{N}_*(\Delta^1)
= \F\cdot\Delta(1)\oplus\F\cdot\Delta(0,1)
\qquad\hbox{and}\qquad\tilde{N}_*(S^1) = \F\cdot\Delta(0,1).$$
The morphism $\tilde{N}_*(S^0)\,\rightarrow\,\tilde{N}_*(\Delta^1)$
takes $\Delta(0)\in\tilde{N}_0(S^0)$ to $\Delta(1)\in\tilde{N}_0(\Delta^1)$.
The morphism $\tilde{N}_*(\Delta^1)\,\rightarrow\,\tilde{N}_*(S^1)$
cancels $\Delta(1)\in\tilde{N}_0(\Delta^1)$
and maps $\Delta(0,1)\in\tilde{N}_1(\Delta^1)$ to $\Delta(0,1)\in\tilde{N}_1(S^1)$.}

\smallskip
We make explicit the cooperations
$\AW(u): \tilde{N}_*(X)\,\rightarrow\,\tilde{N}_*(X)^{\otimes r}$
associated to a fixed surjection $u\in\X(r)_{s-1}$
for $X = S^0,\Delta^1,S^1$.
We assume that this surjection is represented by the sequence
$$u = (u(1),\ldots,u(s),\ldots,u(s+r-1)).$$
The next assertion is immediate:

\paragraph{\sc Fact:}\label{ZeroSphere}
{\it Fix a surjection $u\in\X(r)_0$.
In $\tilde{N}_*(S^0)^{\otimes r}$, we have $\AW(u)(\Delta(0)) = \Delta(0)^{\otimes r}$.
As a consequence, in $\tilde{N}_*(\Delta^1)^{\otimes r}$,
we have $\AW(u)(\Delta(1)) = \Delta(1)^{\otimes r}$.}

\smallskip
Our next purpose is to determine the components of the coproduct
$$\AW(u)(\Delta(0,1)) = \sum\Delta(C_{(1)})\otimes\cdots\otimes\Delta(C_{(r)})$$
in $\tilde{N}_*(\Delta^1)^{\otimes r}$.
In fact, the module $\tilde{N}_*(\Delta^1)^{\otimes r}$ is generated in degree $s$
by the permutations of the tensor $\Delta(0,1)^{\otimes s}\otimes\Delta(1)^{\otimes r-s}$
(where $0\leq s\leq r$).
If $u$ has degree $s-1$, then $\AW(u)(\Delta(0,1))$ has degree $s$.
The component $\Delta(0,1)^{\otimes s}\otimes\Delta(1)^{\otimes r-s}$ of $\AW(u)(\Delta(0,1))$,
where $1\leq s\leq r$, is given by the next lemma.
By $\Sigma_r$-equivariance, this result suffices to determine all components
of the coproduct $\AW(u)(\Delta(0,1))$.

\paragraph{\sc Lemma}\label{OneCone}
{\it Fix a surjection $u\in\X(r)_{s-1}$, where $1\leq s\leq r$.
The coproduct $\AW(u)(\Delta(0,1))\in\tilde{N}_*(\Delta^1)^{\otimes r}$
has no component $\Delta(0,1)^{\otimes s}\otimes\Delta(1)^{\otimes r-s}$
unless the sequences
$$(u(1),u(2),\ldots,u(s))\qquad\hbox{and}\qquad(u(s),\ldots,u(s+r-1))$$
are permutations of $(1,\ldots,s)$ and $(1,\ldots,r)$.
In this case, we have:
$$\AW(u)(\Delta(0,1))
= \sgn(u(1),\ldots,u(s))\cdot\Delta(0,1)^{\otimes s}\otimes\Delta(1)^{\otimes r-s}$$
in $\tilde{N}_*(\Delta^1)^{\otimes r}$.}

\smallskip
As a consequence:

\paragraph{\sc Fact:}\label{OneSphere}
{\it Fix a surjection $u\in\X(r)_{r-1}$.
The coproduct $\AW(u)(\Delta(0,1))$ vanishes in $\tilde{N}_*(S^1)^{\otimes r}$
unless the sequences $(u(1),\ldots,u(r))$ and $(u(r),\ldots,u(r+r-1))$
are both permutations of $(1,\ldots,r)$.
In this case, we have $\AW(u)(\Delta(0,1))
= \sgn(u(1),\ldots,u(r))\cdot\Delta(0,1)^{\otimes r}$.}

\smallskip
The proof of the lemma is postponed to the end of the section.
We now determine the coproducts $w^*: \tilde{N}_*(X)\,\rightarrow\,\tilde{N}_*(X)^{\otimes r}$
associated to an operation $w\in\E(r)_d$.
First, let us record the following fact:

\paragraph{\sc Observation:} {\it Fix $1\leq s\leq r$.
If the surjection $u\in\X(r)_{s-1}$ verifies the condition of the lemma,
then it has the following table arrangement:
$$\left|\matrix{ u(1),\hfill\cr u(2),\hfill\cr\ \vdots\ \hfill\cr u(s-1),\hfill\cr
u(s),\ldots,u(s+r-1).\hfill\cr }\right.$$}

\smallskip
The next assertion follows from this observation.

\paragraph{\sc Claim:} {\it Fix $w = (w_0,\ldots,w_d)\in\E(r)_d$. Assume $d = s-1$.
Let $w'\in\X(r)_d$ be given by the sequence
$$w' = (w_0(1),w_1(1),\ldots,w_{d-1}(1),w_d(1),\ldots,w_d(r)).$$
This surjection $w'$ arises by table reduction from $w$
(we fix the parameters $r_0 = r_1 = \cdots = r_{d-1} = 1$ and $r_d = r$).
The coproduct $\AW(\TR(w)): \tilde{N}_*(X)\,\rightarrow\,\tilde{N}_*(X)^{\otimes r}$,
where $X = \Delta^1$ or $X = S^1$,
reduces to the operations $\AW(w'): \tilde{N}_*(X)\,\rightarrow\,\tilde{N}_*(X)^{\otimes r}$.}

\smallskip
Theorem \ref{SphereModels} follows from this claim, from lemma \ref{OneCone}
and from the facts \ref{ZeroSphere} and \ref{OneSphere} above.
Just observe that the tensor
$${\bf e}_{(1)}\otimes\cdots\otimes{\bf e}_{(s)}
\otimes{\bf c}_{(s+1)}\otimes\cdots\otimes{\bf c}_{(r)}\in C(1)^{\otimes r},$$
where $0\leq s\leq r$,
is dual to $\Delta(0,1)^{\otimes s}\otimes\Delta(1)^{\otimes r-s}$.
To be more precise, because of the commutation rules,
the duality pairing
$$\bigl\langle{\bf e}_{(1)}\otimes\cdots\otimes{\bf e}_{(s)}
\otimes{\bf c}_{(s+1)}\otimes\cdots\otimes{\bf c}_{(r)}\bigl|
\Delta(0,1)^{\otimes s}\otimes\Delta(1)^{\otimes r-s}\bigr\rangle$$
is equal to $(-1)^\sigma$, where $\sigma = s(s-1)/2$.

\paragraph{\it Proof of lemma \ref{OneCone}:}
The surjection is represented by the sequence $(u(1),\ldots,u(s),\ldots,u(s+r-1))$.
We let $\Delta(C_{(1)})\otimes\cdots\otimes\Delta(C_{(r)})$
be the multi-simplex
determined by the interval cuts
$0 = n_0\leq n_1\leq\cdots\leq n_{s+r-1}\leq n_{s+r} = 1$.
To be explicit, we assume $n_i = 0$, for $i = 0,\ldots,j-1$,
and $n_i = 1$, for $i = j,\ldots,s+r$.
By definition, the vertex $\Delta(0)$ is a face
of $\Delta(C_{(i)})$ if $i\in\{\,u(1),\ldots,u(j)\,\}$.
Similarly, the vertex $\Delta(1)$ is a face
of $\Delta(C_{(i)})$ if $i\in\{\,u(j),\ldots,u(s+r-1)\,\}$.
Furthermore, the simplex $\Delta(C_{(i)})$ is degenerate
if the index $s$ occurs twice in $u(1),\ldots,u(j)$ or in $u(j),\ldots,u(s+r-1)$.
Therefore, if the tensor $\Delta(C_{(1)})\otimes\cdots\otimes\Delta(C_{(r)})$
is equal to $\Delta(0,1)^{\otimes s}\otimes\Delta(1)^{\otimes r-s}$
in $\tilde{N}_*(C^1)^{\otimes r}$,
then $(u(1),\ldots,u(j))$ is a permutation of $(1,\ldots,s)$
and $(u(j),\ldots,u(s+r-1))$ is a permutation of $(1,\ldots,r)$.
This property may occur for $j = s$ only
and supposes that the surjection has the following table arrangement:
$$\left|\matrix{ u(1),\hfill\cr\ \vdots\ \hfill\cr
u(s-1),\hfill\cr u(s),\ldots,u(s+r-1).\hfill\cr }\right.$$

Let us now determine the sign
which occurs in the definition of $\AW(u)(\Delta(0,1))$.
The intervals $[n_{i-1},n_i]$ ($i = 1,\ldots,s-1$),
which are reduced to the point $\{0\}$, have length $1$
and have $0$ as a position sign-exponent.
The other intervals are associated to a final occurence of the surjection.
Thus, the interval $[n_{s-1},n_s]$,
which is equal to $[0,1]$, has length 1.
The last intervals $[n_{i-1},n_i]$ (where $i = s+1,\ldots,s+r-1$),
which are reduced to $\{1\}$, have length $0$.
As a consequence, the permutation sign associated to $u$
is identified with the signature $\sgn(u(1),\ldots,u(s))$.
The formula of lemma \ref{OneCone} follows.

\references

\medskip{\parindent=0.5cm\leftskip=0cm

\bibitem{BE}\refto{M. Barratt, P. Eccles},
{\it On $\Gamma_+$-structures. I. A free group functor for stable homotopy theory},
Topology {\bf 13} (1974), 25-45.

\bibitem{Bn}\refto{D.J. Benson},
Representations and cohomology II. Cohomology of groups and modules,
Cambridge Studies in Advanced Mathematics {\bf 31},
Cambridge University Press, 1991.

\bibitem{Bg}\refto{C. Berger},
{\it Op\'erades cellulaires et espaces de lacets it\'er\'es},
Ann. Inst. Fourier {\bf 46} (1996), 1125-1157.

\bibitem{BF}\refto{C. Berger, B. Fresse},
{\it Une d\'ecomposition prismatique de l'op\'erade de Barratt-\break Eccles},
C. R. Acad. Sci. Paris S\'er. I {\bf 335} (2002), 365-370.

\bibitem{BV}\refto{J. Boardman, R. Vogt},
Homotopy invariant algebraic structures on topological spa\-ces,
Lecture Notes in Mathematics {\bf 347},
Springer-Verlag, 1973. 

\bibitem{C}\refto{D. Curtis},
{\it Simplicial homotopy theory},
Advances in Math. {\bf 6} (1971), 107-209.

\bibitem{D}\refto{A. Dold},
{\it\"Uber die Steenrodschen Kohomologieoperationen},
Ann. of Math. {\bf 73} (1961), 258-294.

\bibitem{GZ}\refto{P. Gabriel, M. Zisman},
Calculus of fractions and homotopy theory,
Ergebnisse der Mathematik und ihrer Grenzgebiete {\bf 35},
Springer-Verlag, 1967.

\bibitem{Ge}\refto{E. Getzler},
{\it Cartan homotopy formulas and the Gauss-Manin connection in cyclic homology},
``Quantum deformations of algebras and their representations (Ramat-Gan, 1991/1992; Rehovot, 1991/1992)'',
Israel Math. Conf. Proc. {\bf 7} (1993), 65-78.

\bibitem{GeJ}\refto{E. Getzler, J. Jones},
{\it Operads, homotopy algebra and iterated integrals for double loop spaces},
preprint {\tt arXiv:hep-th/9403055} (1994).

\bibitem{GiK}\refto{V. Ginzburg, M. Kapranov},
{\it Koszul duality for operads},
Duke Math. J. {\bf 76} (1994), 203-272.

\bibitem{GoJ}\refto{P. Goerss, J. Jardine},
Simplicial homotopy theory,
Progress in Mathematics {\bf 174},
Bir\-kh\"au\-ser, 1999.

\bibitem{H}\refto{V. Hinich},
{\it Homological algebra of homotopy algebras},
Comm. Algebra {\bf 25} (1997), 3291-3323.

\bibitem{HS}\refto{V. Hinich, V. Schechtman},
{\it On homotopy limit of homotopy algebras},
in ``$K$-theory, arithmetic and geometry (Moscow, 1984-1986)'',
Lecture Notes in Mathematics {\bf 1289},
Springer-Verlag (1987), 240-264.

\bibitem{Kd}\refto{T. Kadeishvili},
{\it The structure of the $A(\infty)$-algebra, and the Hochschild and Harrison cohomologies},
Trudy Tbiliss. Mat. Inst. Razmadze Akad. Nauk Gruzin. SSR {\bf 91} (1988), 19-27.

\bibitem{K}\refto{T. Kashiwabara},
{\it On the homotopy type of configuration complexes},
{\it in} ``Algebraic topology (Oaxtepec, 1991)'',
Contemp. Math. {\bf 146} (1993), 159-170.

\bibitem{L}\refto{M. Livernet},
Homotopie rationnelle des alg\`ebres sur une op\'erade,
Thesis, Universit\'e Louis Pasteur, Strasbourg, 1998.

\bibitem{ML}\refto{S. Mac Lane},
Homology,
Die Grundlehren der mathematischen Wissenschaften {\bf 114},
Springer-Verlag, 1963.

\bibitem{M}\refto{M. Mandell},
{\it $E_\infty$ algebras and $p$-adic homotopy theory},
Topology {\bf 40} (2001), 43-94.

\bibitem{MSimp}\refto{P. May},
Simplicial objects in algebraic topology,
Van Nostrand, 1967.

\bibitem{MLoop}\reftosame,
The geometry of iterated loop spaces,
Lectures Notes in Mathematics {\bf 271},
Springer-Verlag, 1972.

\bibitem{MCSPrism}\refto{J. McClure, J.H. Smith},
{\it A solution of Deligne's Hochschild cohomology conjecture}
{\it in} ``Recent progress in homotopy theory (Baltimore, 2000)'',
Contemp. Math. {\bf 293}, Amer. Math. Soc. (2002), 153-193.

\bibitem{MCSSeq}\reftosame,
{\it Multivariable cochain operations and little $n$-cubes},
pre\-print {\tt arXiv:math.QA/0106024} (2001).

\bibitem{Q}\refto{D. Quillen},
Homotopical algebra,
Lecture Notes in Mathematics {\bf 43},
Springer-Verlag, 1967.

\bibitem{Sm}\refto{V. Smirnov},
{\it Homotopy theory of coalgebras},
Izv. Akad. Nauk SSSR Ser. Mat. {\bf 49} (1985), 1302-1321.
English translation in Math. USSR-Izv. {\bf 27} (1986), 575-592.

\bibitem{SmCLoop}\reftosame,
{\it On the chain complex of an iterated loop space},
Izv. Akad. Nauk SSSR Ser. Mat. {\bf 53} (1989), 1108-1119.
English translation in Math. USSR-Izv. {\bf 35} (1990), 445-455.

\bibitem{SmHLoop}\reftosame,
{\it The homology of iterated loop spaces},
Forum Math. {\bf 14} (2002), 345-381.

\bibitem{JfS}\refto{J.H. Smith},
{\it Simplicial group models for $\Omega^n\Sigma^n X$},
Israel J. Math. {\bf 66} (1989), 330-350.

\bibitem{JuS}\refto{J.R. Smith},
{\it Operads and algebraic homotopy},
preprint {\tt arXiv:math.AT/0004003}\break (2000).

\bibitem{S}\refto{N. Steenrod},
{\it products of cocycles and extensions of mappings},
Ann. of Math. {\bf 48} (1947), 290-320.

}

\bye